\documentclass[11pt]{article}
\usepackage{amsmath,amssymb,amsthm,mathtools}
\usepackage{dsfont}
\usepackage{aliascnt}
\usepackage{tikz}
\usepackage{quiver}
\usepackage{xcolor}
\usepackage[letterpaper, margin=1in]{geometry}
\usepackage[hidelinks]{hyperref}
\usepackage[nameinlink,noabbrev,capitalize]{cleveref}

\usepackage{graphicx}
\usepackage{float}
\usepackage[section]{placeins}
\usepackage{booktabs}
\usepackage{array}
\usepackage{microtype}
\usepackage{enumitem}
\usepackage{capt-of}

\newtheorem{theorem}{Theorem}[section]

\newaliascnt{claim}{theorem}
\newtheorem{claim}[claim]{Claim}
\aliascntresetthe{claim}

\newaliascnt{lemma}{theorem}
\newtheorem{lemma}[lemma]{Lemma}
\aliascntresetthe{lemma}

\newaliascnt{proposition}{theorem}
\newtheorem{proposition}[proposition]{Proposition}
\aliascntresetthe{proposition}

\newaliascnt{corollary}{theorem}
\newtheorem{corollary}[corollary]{Corollary}
\aliascntresetthe{corollary}

\theoremstyle{definition}
\newaliascnt{definition}{theorem}
\newtheorem{definition}[definition]{Definition}
\aliascntresetthe{definition}

\newaliascnt{example}{theorem}
\newtheorem{example}[example]{Example}
\aliascntresetthe{example}

\theoremstyle{remark}
\newaliascnt{remark}{theorem}
\newtheorem{remark}[remark]{Remark}
\aliascntresetthe{remark}

\theoremstyle{plain}
\newaliascnt{conjecture}{theorem}
\newtheorem{conjecture}[conjecture]{Conjecture}
\aliascntresetthe{conjecture}

\crefname{theorem}{Theorem}{Theorems}
\Crefname{theorem}{Theorem}{Theorems}
\crefname{claim}{Claim}{Claims}
\Crefname{claim}{Claim}{Claims}
\crefname{lemma}{Lemma}{Lemmas}
\Crefname{lemma}{Lemma}{Lemmas}
\crefname{proposition}{Proposition}{Propositions}
\Crefname{proposition}{Proposition}{Propositions}
\crefname{corollary}{Corollary}{Corollaries}
\Crefname{corollary}{Corollary}{Corollaries}
\crefname{definition}{Definition}{Definitions}
\Crefname{definition}{Definition}{Definitions}
\crefname{example}{Example}{Examples}
\Crefname{example}{Example}{Examples}
\crefname{remark}{Remark}{Remarks}
\Crefname{remark}{Remark}{Remarks}
\crefname{conjecture}{Conjecture}{Conjecture}
\Crefname{conjecture}{Conjecture}{Conjectures}
\crefname{section}{Section}{Sections}
\Crefname{section}{Section}{Sections}

\newcommand{\cH}{\mathcal{H}}
\newcommand{\Sym}{\mathrm{Sym}}
\newcommand{\Inv}{\mathrm{Inv}}
\newcommand{\Match}{\mathrm{Match}}
\newcommand{\Z}{\mathbb{Z}}
\newcommand{\N}{\mathbb{N}}
\newcommand{\cL}{\mathcal{L}}
\newcommand{\bF}{\mathbb{F}}
\DeclareMathOperator{\op}{op}

\newcommand{\C}{\mathbb{C}}
\newcommand{\one}{\mathbf{1}}

\newcommand{\mA}{\mathbf{A}}
\newcommand{\mH}{\mathbf{H}}

\newcommand{\mF}{\mathbf{F}}
\newcommand{\mdotG}{\dot{\mathbf{G}}}
\newcommand{\spec}{\operatorname{spec}}
\newcommand{\dH}{d_{\mathrm H}}
\newcommand{\HM}{\cH}

\newcommand{\lamplus}{\lambda_{+}}
\newcommand{\lamminus}{\lambda_{-}}

\newcommand{\Gn}{G_n}

\newcommand{\mdotGn}{\mdotG_n}

\newcommand{\Hone}{\HM_{\mathrm{triv}}}
\newcommand{\Hzero}{\HM_{\mathrm{new}}}
\newcommand{\lambdanew}{\lambda_{\mathrm{new}}}

\newcommand{\ZH}{\mathcal Z_H}
\newcommand{\LambdaH}{\Lambda_H}
\newcommand{\XH}{\mathcal{X}_H}
\newcommand{\mAXH}{\mA_{\XH}}
\newcommand{\mAX}{\mA_X}

\newcommand*{\circled}[1]{\hspace{0.15em}\tikz[baseline=(char.south)]{
    \node[shape=circle,draw,inner sep=0.7pt,scale=0.85] (char) {$\mathsf{#1}$};}\hspace{0.15em}}
\DeclareRobustCommand{\zgzg}{\circled{z}}
\newcommand{\zee}{\mathrm{Z}}

\newcommand{\sm}[1]{\mathcal{#1}}
\newcommand{\gf}[1]{#1}
\newcommand{\E}{\mathbb{E}}
\newcommand{\Prob}{\mathbb{P}}
\newcommand{\pconv}{\overset{\Prob}{\longrightarrow}}
\newcommand{\Tr}{\operatorname{Tr}}
\newcommand{\supp}{\operatorname{supp}}
\newcommand{\diag}{\operatorname{diag}}
\newcommand{\KM}{\mathrm{KM}}
\newcommand{\Tlim}{\mathcal{T}_{\degreeLaw}}
\newcommand{\treed}{\mathcal{T}_{d}}
\newcommand{\nuD}{\nu_{\degreeLaw}}
\newcommand{\rhoDelta}{\rho_{\mindeg}}

\newcommand{\mP}{\mathbf{P}}

\newcommand{\mC}{\mathbf{C}}
\newcommand{\mD}{\mathbf{D}}
\newcommand{\mW}{\mathbf{W}}
\newcommand{\mJ}{\mathbf{J}}
\newcommand{\Mat}{M}
\newcommand{\R}{\mathbb{R}}
\newcommand{\degreeLaw}{D}
\newcommand{\hinv}{h}
\newcommand{\longchain}{L_{\max}}
\newcommand{\pathvalue}{\theta}
\newcommand{\edgefn}{b}
\newcommand{\mWG}{\mW_G}
\newcommand{\mWGn}{\mW_{n}}
\newcommand{\mWT}{\mW_T}
\newcommand{\mDGn}{\mD_{n}}
\newcommand{\mAGn}{\mA_{n}}
\newcommand{\mPGn}{\mP_{n}}
\newcommand{\avgdeg}{\mu}
\newcommand{\mindeg}{\delta}
\newcommand{\maxdeg}{\Delta}
\newcommand{\avgdegn}{\bar d_n}
\newcommand{\avgdegfinite}{\bar d}
\newcommand{\degseq}{\mathbf d^{(n)}}
\newcommand{\portset}{\mathcal J}
\crefname{equation}{Equation}{Equations}
\Crefname{equation}{Equation}{Equations}
\crefname{subsection}{Section}{Sections}
\Crefname{subsection}{Section}{Sections}
\crefname{subsubsection}{Section}{Sections}
\Crefname{subsubsection}{Section}{Sections}
\crefname{figure}{Figure}{Figures}
\Crefname{figure}{Figure}{Figures}
\crefname{table}{Table}{Tables}
\Crefname{table}{Table}{Tables}
\crefname{appendix}{Appendix}{Appendices}
\Crefname{appendix}{Appendix}{Appendices}

\newcommand{\normeig}{\omega}
\graphicspath{{figures/}}
\title{Strong Convergence for Zig-Zag Products and Typical Spectra of Irregular Graphs}
\author{
Alon Dogon\thanks{\texttt{ad4657@columbia.edu}.   Supported by a Junior Fellowship in the Simons Society of Fellows.}\\
{\small Columbia University}
\and
Gal Maor\thanks{\texttt{galmaor@mit.edu}. Supported by the MIT--Kalaniyot Postdoctoral Fellows Program.}\\
{\small MIT}
}

\date{}

\begin{document}
\maketitle
\begin{abstract}
We investigate the following: do existence bounds obtained by interlacing methods describe the spectra of typical random graphs?
Depending on the model, we provide both positive and negative answers.

We establish a Friedman-type theorem for zig-zag products of random regular graphs $G$ with a fixed vertex-transitive gadget $H$, where $G$ is sampled from independent permutations or matchings.
We prove that the new spectrum converges in Hausdorff distance to that of the zig-zag product of the infinite regular tree with $H$. In particular, the new spectral radius converges to the sharp threshold established by Cohen, Cohen, and Maor (FOCS 2024).
We establish the spectral limit using strong convergence and, drawing on quantitative estimates of Chen, Garza-Vargas, Tropp, and van Handel (Annals of Mathematics, 2026), obtain an $O(1/n)$ upper-tail bound in the permutation model.

For the irregular configuration model,
we show that the interlacing bound of Cohen and Maor (RANDOM 2026), asymptotically $2\sqrt{\mu-1}/\mu$ for limiting average degree $\mu$, is atypical: for uniformly bounded degree sequences with a nonconstant limiting degree law of minimum $\delta\ge 2$, a random graph exceeds the interlacing bound by a positive constant with high probability.
This gap is explained by the limiting spectral distribution, whose support edge is $2\sqrt{\delta-1}/\delta$. 
Strikingly, unlike the familiar behavior of random regular graphs, our experiments reveal second eigenvalues that remain far below this limiting edge even for the largest graphs tested. For $\delta\ge3$, we propose that the scarcity of minimum-degree trees explains this behavior; estimates for their size suggest spectral witnesses within $O\left((\log\log n)^{-2}\right)$ of the limiting edge, which is consistent with our numerical observations.
\end{abstract}

\section{Introduction and Main Results}\label{sec:introduction}

The Alon--Boppana bound and Friedman's theorem identify the optimal
asymptotic spectral expansion of regular graphs. For a \(d\)-regular
graph \(G\) on \(n\) vertices, the Alon--Boppana bound~\cite{AB91}
gives
\[
\lambda_2(G)\geq 2\sqrt{d-1}-o_n(1),
\]
where \(\lambda_2(G)\) is the second largest eigenvalue of the
adjacency matrix. The threshold \(2\sqrt{d-1}\) is the spectral
radius of the infinite \(d\)-regular tree, the universal cover of
every connected \(d\)-regular graph. A connected \(d\)-regular graph
is Ramanujan if all its eigenvalues other than the trivial
eigenvalues \(\pm d\) lie in
\([-2\sqrt{d-1},2\sqrt{d-1}]\). Such graphs attain the optimal
threshold without an asymptotic error.

Friedman's theorem~\cite{Fri08,Bor20} shows that almost optimal spectral
expansion is typical: for every fixed \(d\geq3\) and
\(\varepsilon>0\), a random \(d\)-regular graph has all its
eigenvalues other than the top eigenvalue \(d\) in
\[
[-2\sqrt{d-1}-\varepsilon,\,2\sqrt{d-1}+\varepsilon]
\]
with probability tending to one as \(n\to\infty\). Thus the
spectral radius of the regular tree describes both the asymptotic
lower bound for finite graphs and the limiting nontrivial spectral
edge of random regular graphs.

We will distinguish three standard models of random regular (multi-)graphs. The \emph{matchings} model forms a \(d\)-regular graph from \(d\) independent uniform perfect matchings\footnote{We also consider a model we term the \emph{involutions} model, where the perfect matchings are allowed to have self loops.}. For even \(d\), the \emph{permutations} model uses \(d/2\) independent uniform permutations, viewing each permutation edge as undirected (see \cref{subsec:groups_unitaries_permutations}). The \emph{configuration} model instead attaches \(d\) half-edges to every vertex and pairs all half-edges uniformly at random. These models provide closely related ways of sampling regular graphs, but retain different information about how the edges were generated. That distinction becomes important when the sampled graph is used inside another construction.

The zig-zag product, introduced by Reingold, Vadhan and Wigderson~\cite{RVW02}, is a fundamental tool for constructing expanders. Given a \(d\)-regular graph \(G\) and a \(c\)-regular graph \(H\) on \(d\) vertices, the product \(G\zgzg H\) has \(d\cdot|V(G)|\) vertices and degree \(c^2\). The vertices of \(G\zgzg H\) can be viewed as forming a \emph{cloud} of \(H\) for each vertex of \(G\), while an edge in it consists of a step within a cloud (``zig''), a step between clouds determined by \(G\), and another step within a cloud (``zag''). This separates the size of the construction from its degree, and iteration yields explicit constant-degree expanders. For a fixed gadget \(H\), determining the best possible spectral expansion of \(G\zgzg H\) as \(G\) grows is an analogue of the classical Alon-Boppana/Friedman question for regular graphs.
Cohen, Cohen and Maor~\cite{CCM24} established matching lower and existence bounds for this problem, as summarized in the following theorem.
\begin{theorem}[Main theorem of \cite{CCM24}, simplified]
	\label{thm:ccm24_intro}
	Fix a vertex-transitive $c$-regular graph $H$ on $d$ vertices.
    \footnote{Strictly speaking, \cite{CCM24} excludes degenerate cases for which $H^2$ has at least $d/2$ connected components. We omit this mild restriction here for simplicity; in particular, every connected $H$ with $d\ge5$ is admissible.}
	There is an explicitly computable number $\LambdaH$ such that:
	\begin{itemize}
		\item \emph{(``Alon-Boppana'' for zig-zag products.)} For every $d$-regular graph $G$ on $n$ vertices, 
		\begin{equation}\label{eq:ccm_lb}
			\lambda_2(G\zgzg H)\ge \LambdaH-o_n(1).
		\end{equation}
		\item \emph{(Existence of ``zig-zag-Ramanujan'' graphs.)}
		For every size $n$, there exists a $d$-regular graph $G$ in the support of the configuration model such that
		\begin{equation}\label{eq:ccm_ub}
			\lambda_2(G\zgzg H)\le \LambdaH.
		\end{equation}
	\end{itemize}
\end{theorem}
The value \(\LambdaH\) is determined by the (finite) spectrum of \(H\) through a simple calculation. It also has a combinatorial interpretation: let \(\treed\) denote the infinite \(d\)-regular tree and set \(\XH=\treed\zgzg H\). 
Vertex-transitivity of $H$ makes $\XH$ independent, up to isomorphism, of the labeling used to form the product. Moreover, $\XH$ is an infinite cover of $G\zgzg H$ for every connected $d$-regular graph $G$. Although not stated explicitly in \cite{CCM24}, the proof there shows that $\LambdaH$ is the spectral radius of $\XH$. We spell out the details of this identification in \cref{zz:app:cover-spectrum}.

\cref{thm:ccm24_intro} leaves open the ``Friedman-type'' question of how frequently an almost-optimal graph appears. We answer this question in our main result. To state it, we distinguish two invariant subspaces: $\Hone$ is the $d$-dimensional subspace of vectors that assign the same value to corresponding vertices in every cloud, and $\Hzero$ its orthogonal complement\footnote{The terminology is in analogy to the situation of random lifts of graphs, see \cite[Section 5.1.2]{van2025strong}.}.
Let $\zee$ be the adjacency matrix of $G\zgzg H$, when $H$ is fixed and $G$ is picked at random. The restriction $\zee|_{\Hone}$ has deterministic spectrum: while depending on the model, it does not depend on the random realization or the size of $G$. We analyze the spectrum of $\zee|_{\Hzero}$ using the bound $\LambdaH$ developed in \cite{CCM24}.

\begin{theorem}[Informal version of \cref{zz:prop:trivial-block,zz:thm:main,zz:thm:ZZ-spectral-convergence}]\label{thm:zz:intro}
	Fix $H$ as above, and let $\Gn$ be sampled from the $d$-matchings model or, when $d$ is even, the permutations model, with the labeling conventions specified in \cref{sec:main-zigzag}. Let $\zee$ be the adjacency matrix of $\Gn\zgzg H$. Then, as $n$ tends to infinity:
	\begin{enumerate}
		\item We have \label{itm:zgzg_triv} $\spec(\zee|_{\Hone}) = \spec(\mH^2)$ in the matchings model and $\spec(\zee|_{\Hone}) = \spec(\mH\mF\mH)$ in the permutations model, where $\mH$ is the adjacency matrix of $H$ and $\mF$ is the port-flip matrix (defined in \cref{zz:sec:main-results}).
        \item \label{itm:zgzg_friedman}
        \emph{(``Friedman's theorem'' for zig-zag products.)}
        The new spectrum converges in Hausdorff distance to the spectrum of the infinite zig-zag product $\XH$:
        \[
            \dH\left(
                \spec(\zee|_{\Hzero}),
                \spec(\XH)
            \right)\pconv 0.
        \]
        In particular,
        $\left\|\zee|_{\Hzero}\right\|_{\op}\pconv\LambdaH$.
        \item Moreover, in the permutation model, we have for every $\epsilon>0$:
      $$\Prob_{\Gn} \left[\left\|\zee|_{\Hzero}\right\|_{\op}  \geq \LambdaH+ \epsilon \right] = O\left(\epsilon^{-9}\cdot n^{-1}\right).$$
	\end{enumerate}
\end{theorem}

The full details of our zig-zag proof, including a comparison to the configuration model and the results of \cite{CCM24}, appear in \cref{sec:main-zigzag}. At the core of the proof lies \emph{strong convergence}~\cite{van2025strong}, which we briefly survey in \cref{sec:BC}; we show that the zig-zag adjacency matrix can be expressed by a polynomial with fixed matrix coefficients determined by \(H\), an observation also made by O'Donnell and Wu~\cite{OW20}, in their work on derandomizing the Bordenave--Collins theorem.

The formulation of \cref{thm:zz:intro}, unlike \cref{thm:ccm24_intro}, is not framed in terms of $\lambda_2$; this reflects the essential distinction between new and trivial eigenvalues, emphasized in the seminal work of Bordenave and Collins \cite{BC19}.
A function that repeats the same pattern on every cloud is specified by only \(d\) values, independently of \(n\). In each model, the adjacency operator acts on these functions in a way determined only by $H$, the fixed model and edge-labeling convention, contributing \(d\) eigenvalues that can be computed explicitly. The constant function accounts for only one of them, contributing the Perron eigenvalue $\lambda_1=c^2$. The remaining deterministic eigenvalues \emph{can} influence \(\lambda_2\), so removing only the top eigenvalue would not isolate the spectrum controlled by random convergence. Moreover, the deterministic spectrum can differ between the two models even for the same \(H\), as we  demonstrate in \cref{zz:ex:C6-trivial}. 
Together, the deterministic spectrum and the convergence of the largest eigenvalue of $\zee|_{\Hzero}$ determine the limit of $\lambda_2(G\zgzg H)$ in each model. 


While our Friedman-type theorem relies on strong convergence, the existence part of \cref{thm:ccm24_intro} was proved using the methods of \emph{interlacing families} and \emph{finite free convolutions}, developed by Marcus, Spielman, and Srivastava~\cite{MSS2018,marcus2015finite}. These enable establishing spectral existence guarantees, and famously yielded bipartite Ramanujan graphs of every degree and size. In \cite{CCM24}, these were used in combination with the configuration model to analyze the relevant zig-zag operator. Unlike the matching and permutation models considered above, the configuration model extends naturally to irregular degree sequences. Cohen and Maor~\cite{CohenMaor2026} applied the interlacing approach in this setting, establishing an existence bound for every prescribed degree sequence, summarized below.

For a graph $G$, let $\mWG$ denote its symmetric normalized adjacency matrix, whose eigenvalues agree with those of the random-walk matrix. We write $\normeig_2(\mWG)$ for its second largest eigenvalue.
\begin{theorem}[Existence at the average-degree bound {\cite[Theorem~1]{CohenMaor2026}}]
	\label{thm:cohen-maor-existence-intro}
	Let $(d_1,\ldots,d_n)$ be a positive degree sequence with even total degree and average $\avgdegfinite\ge2$. There exists a graph $G$ in the support of the configuration model with these degrees whose second normalized eigenvalue satisfies
	\begin{equation}\label{eq:cohenmaor26}
		\normeig_2(\mWG)\le \frac{2\sqrt{\avgdegfinite-1}}{\avgdegfinite}.
	\end{equation}
\end{theorem}

As we have seen, several graph constructions exhibit a close correspondence between spectral (Alon-Boppana type) lower bounds, interlacing-based existence results, and (Friedman type) typical spectral behavior. 
Examples include regular graphs~\cite{MSS2018} and graph covers~\cite{HPS2018,BC19}; \cref{thm:zz:intro} establishes the corresponding result for zig-zag products.
Would the same hold for the existence result of \cref{thm:cohen-maor-existence-intro}? As observed in \cite{CohenMaor2026}, the expression in \cref{eq:cohenmaor26} is \emph{not} a universal Alon--Boppana lower bound (albeit the visual similarity), since for some degree sequences better expansion is possible.

To state our results fully, we introduce some notations. For each $n$, fix a positive degree sequence $\degseq=(d_1,\ldots,d_n)$ with even total degree and average $\avgdegfinite$. Assume that the degrees are uniformly bounded, independently of $n$, and that their empirical distribution converges to a law $\degreeLaw$. Let $\Gn$ be sampled from the configuration model with these degrees, and let $\mWGn:=\mW_{\Gn}$. Write
\[
\avgdeg:=\E\degreeLaw,
\qquad
\mindeg:=\min\{d:\Prob(\degreeLaw=d)>0\},
\]
and assume throughout that $\mindeg\ge2$.

We refute the analogous Friedman-type expectation for irregular configuration graphs. Whenever the limiting degree law $\degreeLaw$ is nonconstant, the second normalized eigenvalue of a random configuration graph exceeds the Cohen--Maor existence bound by a fixed positive amount with probability tending to one. Thus, the familiar agreement between existence and typicality, seen above for regular graphs, graph covers, and zig-zag products, breaks down in the irregular setting.

The relevant limiting object is the \emph{unimodular Galton--Watson tree} $\Tlim$ (see, e.g.,~\cite[Theorem~4.1]{Hofstad2024}). Its root $r$ has degree distributed according to $\degreeLaw$, while each non-root vertex independently has $\degreeLaw^*-1$ children, where $\degreeLaw^*$ is the size-biased degree law
\[
\Prob(\degreeLaw^*=c)=\frac{c\Prob(\degreeLaw=c)}{\avgdeg}.
\]
This size bias reflects the greater probability of reaching a higher-degree vertex when exploring along an edge.
We denote by $p_k(\Tlim,r)$ the probability that simple random walk on $\Tlim$, started at $r$, is back at $r$ after exactly $k$ steps.
The following theorem summarizes the limiting spectral law and its implications for typical spectral expansion. Various parts of it and its building blocks were already known, and can be attributed to \cite{Hofstad2024}, \cite{BordenaveLelarge2010} (for convergence) and \cite{MuellerWiegel2020} (for the spectral edge).

\begin{theorem}[Limiting spectral law and atypicality of the average-degree bound]
    \label{thm:irreg_weak_intro}
    The empirical spectral measure of $\mWGn$ converges weakly in probability to a measure $\nuD$, where $\nuD$ is the unique probability measure on $[-1,1]$ whose moments are
    \begin{equation}
        m_k(\nuD)=\E\left[p_k(\Tlim,r)\right].
    \end{equation}
    The limiting law $\nuD$ is symmetric and has support edge
    \begin{equation}\label{eq:irreg_bound_intro}
        \rhoDelta:=\frac{2\sqrt{\mindeg-1}}{\mindeg}.
    \end{equation}
    Moreover, and in contrast to the existence guarantee of \cref{thm:cohen-maor-existence-intro}, if $\degreeLaw$ is non-constant, there exists $\varepsilon_{\degreeLaw}>0$ such that
    \[
        \Prob\left[
            \normeig_2(\mWGn)
            \ge
            \frac{2\sqrt{\avgdegfinite-1}}{\avgdegfinite}
            +\varepsilon_{\degreeLaw}
        \right]
        \underset{n\to\infty}{\longrightarrow}1.
    \]
\end{theorem}

For $\mindeg\ge3$, this does not establish convergence of $\normeig_2(\mWGn)$ to $\rhoDelta$. In fact, \cref{cfg:cor:atypical} gives a stronger threshold statement: for every fixed $t<\rhoDelta$, a positive proportion of the eigenvalues of $\mWGn$ exceed $t$ with probability tending to one.

The existence bound in \cref{thm:cohen-maor-existence-intro} is not optimal for certain degree sequences, as noted in \cite{CohenMaor2026}; our result shows that it is nevertheless atypically good for random pairings when $\degreeLaw$ is nonconstant. Thus the expression in \cref{eq:cohenmaor26} serves neither as a universal Alon--Boppana threshold nor as an asymptotically typical (``Friedman'') value.
As we discuss in \cref{cfg:subsec:rotation-free-benchmark}, the free-probability calculation underlying \cite{CohenMaor2026} is inspired by free projections and convolution semigroups (see \cite[Lecture~14]{NS06}). In the limit $\avgdegfinite\to\avgdeg$, it produces the normalized Kesten--McKay law $\KM_\avgdeg$, a benchmark depending only on the average degree rather than the full degree distribution.

However, this free-probability benchmark does not, in general, describe the typical limiting spectral distribution of random configuration graphs. Indeed, while their second moments agree, we show that
\[
m_4(\nuD)-m_4(\KM_\avgdeg)
=\frac1\avgdeg\left(
\E\left[\frac1\degreeLaw\right]-\frac1\avgdeg
\right)^2,
\]
which is strictly positive whenever $\degreeLaw$ is nonconstant. This provides a simple, explicit witness to the limitations of the free-probability benchmark for typical random configuration graphs.

The limiting support edge raises a finite-size puzzle. Consider the degree law with $\Prob(\degreeLaw=3)=2/3$ and $\Prob(\degreeLaw=12)=1/3$, so that $\avgdeg=6$ and $\mindeg=3$. The Cohen--Maor existence bound is then $2\sqrt{5}/6\approx0.745$, whereas the limiting support edge is $\rhoDelta=2\sqrt{2}/3\approx0.943$. Yet in our experiments, the median $\normeig_2$ increases only from $0.765$ at $n=500$ to $0.809$ at $n=50{,}000$; see \cref{cfg:fig:min3-convergence}. We propose a concrete explanation: sufficiently large regions consisting entirely of minimum-degree vertices may be too scarce in finite graphs to produce eigenvalues near $\rhoDelta$.

In \cref{cfg:sec:experiments}, we study these regions as \emph{spectral witnesses}: local subgraphs supporting test vectors with large Rayleigh quotients. When $\mindeg=2$, these are long degree-$2$ chains, whose maximum length is heuristically of order $\log n$. When $\mindeg\ge3$, they are branching trees, whose maximum depth grows only like $\log\log n$. The corresponding Rayleigh quotients suggest gaps from $\rhoDelta$ of order $(\log n)^{-2}$ and $(\log\log n)^{-2}$, respectively. These heuristic scales provide a possible explanation for the strikingly slow finite-size behavior.

\section{Strong convergence and its consequences}\label{sec:BC}

Strong convergence has proved useful in studying spectral properties of random graph lifts \cite{BC19} and spectral gaps of hyperbolic surfaces \cite{hide2023near}. In this section, we recall the basic definitions and results needed for our purposes, following the surveys \cite{van2025strong,magee2025strong}.

\subsection{Groups, unitaries, permutations and graphs}\label{subsec:groups_unitaries_permutations}

We will work in the general setting of a finitely generated group $\Gamma$ with a specified  finite symmetric generating set $S \subset \Gamma$.
This will allow us to provide uniform arguments for both the random permutation and matching models.

For the permutation model, the relevant group is a \emph{free group}: For $r \in \N$, we let $\bF_r = \langle s_1, \dots, s_r \rangle$ be the free group generated by a set of formal generators $s_1, \dots, s_r$.
For convenience, we let $s_{2r - i+1} = s_i^{-1}$ for $1 \leq i \leq r$, and let $S = \{ s_1, \dots, s_{2r} \}$ be a finite \emph{symmetric} generating set for $\bF_r$.
We also denote $\bar{i} = 2r-i + 1$ for $1 \leq i \leq 2r$, so that $s_{\bar{i}} = s_i^{-1}$.
Actions of the free group yield $2r$-regular labeled graphs, as explained later on.

 For the matchings and involutions model, the relevant group is the \textit{free product}:
$$ \Z_2^{*r} = \underbrace{\Z_2 * \Z_2 * \cdots * \Z_2}_{r\text{ times}} = \langle s_1, \dots , s_r\rangle.$$
Here each $s_i$ is an involution. We naturally view $\Z_2^{*r}$ as a quotient group of $\bF_r$ by mapping a generator and its inverse to the corresponding involution.
Setting $S = \{s_1, \dots, s_r \}$ we get a symmetric generating set as $s_i^{-1}  = s_i$ in this case.
Later on, we will see that actions of $\Z_2^{*r}$ yield $r$-regular labeled graphs.
Except for the difference in regularity, all that follows will apply equally well to both.

 A \textbf{\textit{unitary $r$-tuple}} is a tuple $\boldsymbol{u}=(u_1, \dots, u_r) \in U(\cH)^r$ of unitaries on some (possibly infinite dimensional) Hilbert-space $\cH$.
 We will freely and implicitly use the bijection between a unitary $r$-tuple $\boldsymbol{u} = (u_1, \dots, u_r) \in U(\cH)^r$ and unitary representations $\pi_{\boldsymbol{u}}: \bF_r \to U(\cH)$:
Given an $r$-tuple $\boldsymbol{u}$ and a word $w \in \bF_r$, $\pi_{\boldsymbol{u}}(w)$ is given by substituting $s_i$ with $u_i$ in $w$, conversely, given $\pi: \bF_r \to U(\cH)$, we recover a unitary $r$-tuple by letting $u_i = \pi(s_i)$.
Note that the representation $\pi_{\boldsymbol{u}}$ factors to a representation $\pi_{\boldsymbol{u}}: \Z_2^{*r} \to U(\cH)$ if and only if $u_i^2 = 1$ for all $i$, we call such tuples \emph{\textbf{involutory}}.

By a \textbf{\textit{random unitary $r$-tuple}} in $U(n)$, we mean a random variable with values in $U(n)^r$, equivalently, a probability measure on $U(n)^r$.
We will denote such a random tuple by a bold symbol $\boldsymbol{u}^n=(u_1^n,\dots, u_r^n)$.
By a \textbf{\textit{random permutational $r$-tuple}}, we mean a random variable with values in $\Sym(n)^r$.
We will view $\Sym(n) \subset U(n)$ naturally via permutation matrices, and so a permutational $r$-tuple naturally gives rise to a unitary $r$-tuple.
Denote by $ \Inv(n) = \{ \sigma | \sigma^2 = 1 \}\subset \Sym(n)$ the set of \textbf{involutions}. An involution is a \textbf{matching} if it is fixed point-free, denote them by $\Match(n)$.
Using Schreier graphs recalled below, $r$-permutations give rise to $2r$-regular graphs, while $r$-involutions and $r$-random matchings give rise to $r$-regular graphs, with the latter being self-loop free.

\begin{definition}[Schreier graphs]
Let $\Gamma$ be a finitely generated group with a specified symmetric generating set $S$.
Let $\pi: \Gamma \to \Sym(V)$ be a group action, with $V$ a set. 
The associated \textbf{\textit{Schreier graph}} is the $S$-labeled graph $G_\pi$ with vertex set $V$, where two vertices $x, y \in V$ are connected by an edge labeled $s \in S$ if $\pi(s)\cdot x = y$.

While the Schreier graph can be considered as a directed graph, we identify opposite directed edges $x \xrightarrow[]{s} y$ and $x \xleftarrow[]{s^{-1}}y$ to a single undirected edge.
\end{definition}

In the case $\Gamma = \bF_r$, given a permutational $r$-tuple $\boldsymbol{\sigma} = (\sigma_1,\dots,\sigma_r)\in\Sym(n)^r$, the associated Schreier graph $G_{\boldsymbol{\sigma}}$ is a $2r$-regular labeled graph given by connecting two vertices $x, y \in [n]$ by an edge labeled $1 \leq i \leq 2r$ if $\sigma_i\cdot x = y$.
In the case $\Gamma = \Z_2^{*r}$, given an involutory $r$-tuple $\boldsymbol{\sigma} = (\sigma_1,\dots,\sigma_r)\in\Inv(n)^r$, the associated Schreier graph is an $r$-regular labeled graph given by connecting two vertices $x, y \in [n]$ by an edge labeled $1 \leq i \leq r$ if $\sigma_i\cdot x = y$.
Note that $\sigma_i$ are all matchings if and only if $G_{\boldsymbol{\sigma}}$ is self-loop free. 

We now express the permutation and matching models introduced in the introduction as random group actions, together with the closely related involution model.

\begin{definition}[Random actions and random graph models]
For each $n$, consider the following random permutational $r$-tuples and their corresponding actions.
\begin{enumerate}
    \item[\textbf{(P)}]
    \textbf{Independent permutations model}
Here we simply pick $r$-random permutations $\sigma_1^n, \dots \sigma_r^n \in \Sym(n)$ independently, uniformly at random. We then let $u_i^n$ be the permutation matrix associated with $\sigma_i^n$.
This is the same as picking a uniformly random homomorphism $\pi: \bF_r \to \Sym(n)$.
We let $G_{\boldsymbol{\sigma}^n}$ be the resulting random $2r$-regular graph.

    \item[\textbf{(M)}]
    \textbf{Independent matchings model}
    Here, assuming $n$ is even, we pick  $\tau_1^n, \dots, \tau^n_r \in \Match(n)$ independently and let $u_i^n$ be the corresponding self adjoint permutation matrices.
    Note that the random matchings model is the same as picking a uniformly random homomorphism $\pi: \Z_2^{*r} \to \Sym(n)$, for which the standard generators of $\Z_2^{*r}$ act without fixed-points.
    We let $G_{\boldsymbol{\tau}^n}$ be the resulting random $r$-regular graph.
    
    \item[\textbf{(I)}]
    \textbf{Independent involutions model}
    The same as the independent matching model, without the fixed point free condition. 
\end{enumerate}
\end{definition}

Finally, in order to define strong convergence, we recall the relevant limiting objects.
For a group $\Gamma$, we consider the Hilbert-space $\ell^2(\Gamma)$ of square summable complex valued functions on $\Gamma$.
The \textbf{left regular representation} of $\Gamma$ is the unitary representation $\lambda: \Gamma \to U(\ell^2 (\Gamma))$ given by left multiplication.
The \textbf{free Haar unitaries} are the unitaries $u_1, \dots, u_r \in U(\ell^2(\bF_r))$ given by $u_i = \lambda(s_i)$.
Similarly, the \textbf{free Haar involutions} are $v_1, \dots, v_r \in U(\ell^2(\Z_2^{*r})) $ given by $v_i = \lambda(s_i)$.

\subsection{Strong convergence}

We denote by $\C\langle\boldsymbol{x}, \boldsymbol{x}^* \rangle = \C\langle x_1, \dots ,x_r, x_1^*, \dots, x_r^* \rangle$ the $*$-algebra of complex $*$-polynomials with $r$ non-commuting variables. 
Given a unitary $r$-tuple $\boldsymbol{u}$ on $\cH$ and $P \in   \C\langle\boldsymbol{x}, \boldsymbol{x}^* \rangle$ we will write $P(\boldsymbol{u}) = P(\boldsymbol{u}, \boldsymbol{u}^*) \in \cL(\cH)$ for the resulting operator obtained from plugging in the unitaries $u_1, \dots, u_r$ and their inverses into $P$.
The following fundamental definition is due to Haagerup and Thorbj{\o}rnsen \cite{haagerup2005new}.

 \begin{definition}\label{def:storng_conv}
 Fix $r \in \N$ and let $\cH$ be a Hilbert-space. A sequence of random unitary $r$-tuples $\boldsymbol{u}^n = (u_1^n, \dots, u_r^n)$ in $U(n)$ of increasing dimension is said to \emph{strongly converge} to a unitary $r$-tuple  $\boldsymbol{u} = (u_1, \dots ,u_r)$ in $U(\cH)$ if  for every $*$-polynomial $P \in \C\langle \boldsymbol{x}, \boldsymbol{x}^* \rangle$, the following converges in probability
 $$ \| P(\boldsymbol{u}^n) \|_{\op} \underset{n\to\infty}{\pconv} \| P(\boldsymbol{u}) \|_{\op}.$$

In the setting of a general group $\Gamma$ with a finite generating set $S$, a sequence of unitary representations $\pi_n: \Gamma \to U(n)$ converges strongly to a unitary representation $\pi:\Gamma \to U(\cH)$ if after enumerating the generators $S = \{ s_1, \dots ,s_r\}$, the corresponding unitary tuples $u_i^n = \pi_n(s_i)$ converge strongly to the unitary tuple $u_i = \pi(s_i)$.
 \end{definition}

By breakthrough work of \cite{BC19}, alongside the later influential work \cite{chennew} and subsequent work \cite{barbieri2026strong}, we have the following:

\begin{theorem}[\cite{BC19, chennew, barbieri2026strong}]
\label{thm:strong_conv_qualitative}
Let $\boldsymbol{u}^n = (u_1^n, \dots, u^n_r)$ be a random matrix tuple obtained from either the random permutation, matching or involution models. Let $\Gamma$ be either $\bF_r$ in the permutation case, or $\Z_2^{*r}$ in the matching and involution case. 
Then the tuples
$$ \boldsymbol{u}^n|_{\mathds 1^\perp} = (u_1^n|_{\mathds 1^\perp}, \dots, u_r^n|_{\mathds 1^\perp})$$
converge strongly to the regular representation of $\Gamma$.
\end{theorem}

 It is known that strong convergence can automatically be upgraded to include matrix valued polynomials, as well as control on the Hausdorff distance $\dH$ to the full spectrum of the limiting operator.
 Recall that for $A,B \subset \R$, their Hausdorff distance is: 
$$ \dH(A,B) = \max\{ \sup_{x \in A} (d(x, B)), \sup_{y \in B}(d(y, A))\},$$
where $d(x, B) = \inf_{y \in B}|x - y|$ is the distance function.

 \begin{proposition}[Matrix valued strong convergence~\cite{van2025strong}, Lemmas 2.11, 2.16] \label{lem:matrix_valued_conv}
  A sequence of $r$-tuples of unitary matrices $\boldsymbol{u}^n$ converges strongly to $\boldsymbol{u}$ if and only if for every $d \in \N$ and every (possibly non-self adjoint) matrix valued polynomial $P \in \Mat_d(\C) \otimes \C\langle \boldsymbol{x}, \boldsymbol{x}^* \rangle$, the following converges in probability:
  $$ \| P(\boldsymbol{u}^n) \|_{\op} \underset{n\to\infty}{\pconv} \| P(\boldsymbol{u}) \|_{\op}. $$
  If $P(\boldsymbol{u}^n)$ and $P(\boldsymbol{u})$ are self-adjoint for all $n$, then
  $$\dH\left(\spec(P(\boldsymbol{u}^n)), \spec(P(\boldsymbol{u}))\right) \underset{n\to\infty}{\pconv} 0.$$
 \end{proposition}

Part of the novelty in the work \cite{chennew} is yielding strong effective bounds on the convergence rate in  \cref{thm:strong_conv_qualitative}, in the case of the permutations model.
 For a matrix polynomial $P \in \Mat_D(\C) \otimes  \C\langle\boldsymbol{x}, \boldsymbol{x}^* \rangle$, let $\| P \|_1 = \sum_{w} \|A_w \|_{\op}$ where $P=  \sum_{w} A_w \otimes w(\boldsymbol{x}, \boldsymbol{x}^*)$ is written uniquely as a sum of non-commutative monomials.
We have that $\| \pi(P) \|_{\op} \leq \| P \|_1 $ for any unitary representation $\pi: \bF_r \to U(\cH)$, so $\| P \|_{\Mat_D(\C) \otimes C^*(\bF_r)} \leq \| P \|_1$ in the terminology of \cite[Section 3.3]{chennew}.
 
\begin{theorem}[Simplified form of {\cite[Theorem~3.9]{chennew}}]
\label{thm:strong_conv_effective}
Fix $r\ge2$. For any $P \in \Mat_D(\C) \otimes  \C\langle\boldsymbol{x}, \boldsymbol{x}^* \rangle$ a self adjoint matrix valued polynomial of degree $q$. Set $K = \|P\|_1$, then for all $n$ and $0<\epsilon<K$:

  \begin{equation*}
      \Prob_{\boldsymbol{u^n}} \left[ \| P(\boldsymbol{u}^n|_{\mathds{1}^\perp}) \| \geq \|P(\boldsymbol{u}) \| + \epsilon \right] \lesssim \frac{D}{n} \left( \frac{Kq\log(r)}{\epsilon}\right)^8\log\left(\frac{eK}{\epsilon}\right)
  \end{equation*}

  Where $\boldsymbol{u}$ is the left regular representation of $\bF_r$, and $\boldsymbol{u}^n$ are sampled from the permutation model.
\end{theorem}

\section{Strong Convergence of Zig-Zag Products}\label{sec:main-zigzag}

 In this section, we give the definition of zig-zag products for general group actions using the terminology established in \cref{sec:BC}.
 This will help capture both the permutational and matching models under one framework, and recovers the original construction of \cite{RVW02}.
 
Throughout, we fix a finitely generated group $\Gamma$ with a finite symmetric generating set $S$.
For $s_1, s_2 \in S$, we consider the elementary matrix $E_{s_1,s_2} \in \cL(\ell^2(S))$ given by taking the basis vector $\mathds{1}_{s_1}$ to $\mathds{1}_{s_2}$.

\begin{definition}[Rotation map]
    Let $\pi: \Gamma \to \Sym(V)$ be a group action, $V$ a (possibly infinite) set.
    The  \emph{\textbf{rotation map}} of the Schreier graph $G_\pi$ is the permutation $\mdotG_{\boldsymbol{\pi}}: S \times V \to S \times V$ given by $\mdotG_{\boldsymbol{\pi}}(s, x) = (s^{-1}, \pi(s)\cdot x)$.
\end{definition}

It is easy to see that the rotation map is in fact an involution.
Note that if $\Gamma =\bF_r$, $\sigma_1, \dots ,\sigma_r \in \Sym(V)$ describe the action, then $\mdotG_\pi (i,x) = (\bar i, \sigma_i\cdot x)$.
This recovers the definition of rotation maps used in \cite[Definition 2.1]{RVW02}: 
Indeed, the edge $x \xrightarrow[]{i} \sigma_ix$ is also the edge $x \xleftarrow[]{\bar i} \sigma_i x$ in the Schreier graph, as we identified opposite edges.
In the case the graph is described by $r$-matchings $\tau_1, \dots, \tau_r \in \Match(V)$ as a Schreier graph of $\Gamma = \Z_2^{*r}$, then  $\mdotG_\pi (i,x) = (i, \tau_i\cdot x)$, again recovering \cite{RVW02}.

We are now in position to recall the definition of the zig-zag product of \cite{RVW02}.
Let $H$ be a $c$-regular graph whose vertex set is $S$, and denote its adjacency operator by
\[
\mA_H\in \cL(\ell^2(S)).
\]
Thus the vertices of $H$ are identified with the labels, or ports, of the Schreier graphs associated with actions of $\Gamma$.

\begin{definition}[Zig-Zag product]\label{zz:def:action-zigzag}
Let $\Gamma$ be a finitely generated group with a symmetric generating set $S$.
Let $H$ be a $c$-regular graph whose vertex set is $S$, for some $c \in \N$.

Given an action $\pi:\Gamma\to\Sym(V)$ on a set $V$, and let $G_\pi$ be its Schreier graph.  The \emph{zig-zag product of the action $\pi$ with the gadget $H$}, denoted  $G_\pi\zgzg H$, is the graph on $S\times V$ whose edges are given by such paths:
\[\begin{tikzcd}
	{(x,\textcolor{blue}{s_1})} && {(x,\textcolor{blue}{s_2}) } \\
	\\
	{(s_2 \cdot x, \textcolor{blue}{s_3})} && {(s_2\cdot x, \textcolor{blue}{s_2^{-1}})}
	\arrow["{zig}", color={blue}, from=1-1, to=1-3]
	\arrow["zig-zag", color={red}, from=1-1, to=3-1]
	\arrow["{ \mdotG_{\pi}}", from=1-3, to=3-3]
	\arrow["{zag}", color={blue}, from=3-3, to=3-1]
\end{tikzcd}\]

Where $s_1, s_2, s_3 \in S$ are such that $s_1, s_2$ are connected and $s_2^{-1}$ is connected to $s_3$ in $H$.
\end{definition}

Thus a zig-zag step consists of a step inside $H$, followed by a Schreier-graph step and the corresponding inversion of the port, followed by another step inside $H$. 
The zig-zag product is thus a $c^2$-regular graph.
Again, by opening up the definition in the case of $\Gamma = \bF_r$ and a Schreier graph given by $r$-permutations, the definition above fits the original definition of a zig-zag product (see \cite[Definition 4.4.]{CCM24}).
The same is true for $\Gamma= \Z_2 ^{*r}$ and $r$-matchings.

\begin{remark}
    While zig-zag products arising from Schreier graphs of group actions are all zig-zag products in the sense of \cite{RVW02}, the converse is not true.
    Generally, zig-zag products can be constructed from any $d$-regular graph $G=(V,E)$ for which every vertex $v$ labels its neighbors by $v[i], 1 \leq i\leq d$.
    In Schreier graphs, we labeled the neighbors of $v$ by the generating set $S \leq \Gamma$, and in such a way that $v= (v[s])[s^{-1}]$.
    Such a consistency condition need not hold for generally labeled graphs as in \cite{RVW02}.
    For example, this is the case in the configuration model.
    The local consistency present in our labelings is crucial for applying existing strong convergence results to control the spectrum of the zig-zag product. 
\end{remark}

The adjacency operator of the zig-zag product admits a simple formula in terms of the rotation map (see \cite[Section 4.1.2]{CCM24}).
Under the isomorphism $\ell^2(S\times V) = \ell^2(S) \otimes \ell^2(V)$ as Hilbert spaces, we have:
\begin{equation}
\mA_{G_\pi\zgzg H}
:=
(\mA_H\otimes I_V)\,
{\mdotG}_\pi\,
(\mA_H\otimes I_V).
\label{zz:eq:zz-matrix}
\end{equation}
The spectrum of a zig-zag product $G_\pi \zgzg H$ will always contain eigenvalues coming from $H$, which we compute in \cref{zz:sec:trivial-spectrum}. We call these the \emph{trivial eigenvalues}, in analogy to the case of lifts of graphs (see \cite[Section 5.1]{van2025strong}).

\begin{definition}[Trivial and new spectra]
Consider the  following subspaces invariant under $\mA_{G_{\pi} \zgzg H}$.
\[
{\Hone}
:=
\ell^2(S)\otimes\operatorname{span}\{\mathds{1}\},
\qquad
{\Hzero}
:=
\ell^2(S)\otimes{\mathds{1}}^\perp.
\]
The spectrum of the operator $\mA _{G_{\pi} \zgzg H}$ restricted to $\Hone$ is called the \textbf{\textit{trivial spectrum}} and is denoted $\spec_{\textrm{triv}}(G_\pi \zgzg H)$.
The spectrum of the operator $\mA _{G_{\pi} \zgzg H}$ restricted to $\Hzero$ is called the \textbf{\textit{new spectrum}} and is denoted $\spec_{\textrm{new}}(G_\pi \zgzg H)$.
We define the \emph{new spectral radius} $\lambdanew(G_\pi \zgzg H)$ by the max absolute value of $\lambda\in\spec_{\textrm{new}}(G_\pi\zgzg H)$. Equivalently, 
$
\lambdanew(G_\pi \zgzg H)
:=
\left\|\mA_{G_\pi\zgzg H}|_{\Hzero}\right\|_{\op}.
$

\end{definition}

\subsection{The main results}\label{zz:sec:main-results}

Below is our main result mentioned in the introduction spelled out in full.
We now specialize to the permutation, matching and involution models, to derive the explicit spectral edges.
For the explicit spectral edge below, we impose the same standing hypotheses as in \cite[Theorem~5.2]{CCM24}, and exclude the degenerate vertex-transitive graphs $H$ isolated in \cite[Definition~5.1]{CCM24}.

Given $H$ $c$-regular on $d$-vertices which is \emph{good} in the sense of \cite{CCM24}, set
\[
h_H(x):=\det(xI_d-\mA_H^2),
\]
and, for $x>c^2$, define
\begin{equation}
\ZH(x)
:=
x\sqrt{
1-\frac{d\,h_H(x)}
{x\,h_H'(x)}
},
\qquad
\LambdaH
:=
\min_{x>c^2}\ZH(x).
\label{zz:eq:ZH-definition}
\end{equation}

The formula originates in \cite[Section~5.1]{CCM24} from a generating-function analysis of closed walks on the infinite zig-zag product $\treed\zgzg H$. The matching existence bound is established using finite free multiplicative convolution \cite[Section~5.3]{CCM24}. We use unnormalized adjacency throughout, so $\LambdaH$ is $c^2$ times the corresponding normalized quantity in the zig-zag analysis of \cite{CCM24}.

For the symmetric generating set $S$, define the \emph{port-flip involution} $\mF_S\in\cL(\ell^2(S))$ by
\[
    \mF_S := \sum_{s\in S} E_{s,s^{-1}},
\]
thus $\mF_S\mathds{1}_s=\mathds{1}_{s^{-1}}$. In the permutation model, it exchanges port $i$ with port $\bar i$; in the matching and involution models, $\mF_S=I$. The first theorem computes the trivial spectrum and the typical edge of the new spectrum for random zig-zag products, with a quantitative bound in the permutation model.

\begin{theorem}[Trivial spectrum and upper-tail bounds]\label{zz:thm:main}
    Let $H$ be a $c$-regular vertex-transitive graph with $2r$-vertices in the first clause, and $r$ vertices in the second clause, which is good in the sense of \cite{CCM24}.
    
    \begin{enumerate}
        \item  For $\Gn$ sampled in the permutation model, we have $\spec_{\textrm{triv}}(\Gn \zgzg H) = \spec(\mA_H\mF_S \mA_H)$ and for every $\epsilon>0$:
      $$\Prob_{\Gn} \left[  \lambdanew(\Gn \zgzg H)  \geq \LambdaH+ \epsilon \right] = O(1/(n\cdot \epsilon^{9})) $$
        Where the constants only depend on $H$.
    \item For $\Gn$ sampled in the matching or involution model, we have $\spec_{\textrm{triv}}(\Gn \zgzg H) = \spec(\mA_H^2)$ and  the edge of the new spectrum satisfies for every $\epsilon > 0$:
    $$\Prob_{\Gn} \left[  \lambdanew(\Gn \zgzg H)  \geq \LambdaH+ \epsilon \right] = o(1) $$
    \end{enumerate}
    \end{theorem}

A key point in the above theorem is that $\LambdaH$ turns out to be the spectral radius of the zig-zag product with a tree $\treed\zgzg H$, which is an infinite $c^2$-regular graph. The second main theorem controls the \emph{Hausdorff distance} to the spectrum of $\treed\zgzg H$, thus controlling the full spectrum of random zig-zag products.
\begin{theorem}[Full spectral convergence]\label{zz:thm:main_full_spec}
    Let $H$ be as in \cref{zz:thm:main}. Then
    \begin{enumerate}
        \item For $\Gn$ sampled in the permutation model, the following converges in probability
$$ \dH( \spec(\mA_{\Gn \zgzg H}),  \spec(\mA_H\mF_S\mA_H)  \cup \spec(\mA_{\mathcal T_{2r} \zgzg H})) \underset{n\to\infty}{\pconv} 0$$
        \item For $\Gn$ sampled in the matching and involution model, the following converges in probability
$$ \dH( \spec(\mA_{\Gn \zgzg H}),  \spec(\mA_H^2)  \cup \spec(\mA_{\mathcal T_{r} \zgzg H})) \underset{n\to\infty}{\pconv} 0$$
    \end{enumerate}
\end{theorem}

Although \cref{zz:thm:main_full_spec} controls the full spectrum of random zig-zag products in terms of the spectrum of $\treed\zgzg H$, we do not explicitly compute $\spec(\mA_{\treed\zgzg H})$.
It is a consequence of \cite{CCM24} that we have an explicit formula for its edge as in \cref{zz:thm:main}.
In the following section we prove these results.

\subsection{The new spectrum of Random Zig-Zag Products}
\label{zz:sec:spectral-convergence}
Throughout the section, let $\Gamma$ be a finitely generated group with a fixed finite symmetric generating set $S$ and $H$ be a $c$-regular gadget graph on the generators $S$.

In order to control the new spectrum of typical zig-zag product $G_\pi \zgzg H$, we will use strong convergence and identify the limit object.

\begin{definition}
    Given $\Gamma, S, H$ as above, the (left) \textbf{\emph{regular zig-zag-product}} is the zig-zag product $G_\lambda \zgzg H$ obtained from the left-regular action $\lambda: \Gamma \to \Sym(\Gamma)$ given by $\lambda(g)\cdot h = gh$. 
\end{definition}

When $\Gamma$ is infinite, the regular zig-zag product is an infinite $c^2$-regular graph.
In the case of $\Gamma = \bF_r$ (respectively $\Gamma = \Z_2^{*r}$), the regular zig-zag product $\mathcal T_{2r} \zgzg H$ (respectively $\mathcal T_r \zgzg H$) plays the same role as the infinite $2r$-regular (respectively $r$-regular) tree does for the family of all $2r$-regular (respectively $r$-regular) graphs.
Indeed, if $\pi:\Gamma\to\Sym(V)$ is transitive, the natural port-preserving covering map from the regular Schreier graph of $\Gamma$ to $G_\pi$ induces a covering map from $G_\lambda\zgzg H$ to $G_\pi\zgzg H$.
This makes it the natural candidate for a limiting object of such zig-zag products.

The following shows that under the assumption of strong convergence, the \emph{entire} spectrum of the regular zig-zag product controls in Hausdorff distance the \textit{new spectrum} of random zig-zag products with fixed $\Gamma,S,H$, yielding the following "Friedman-type" theorem.

\begin{theorem}[Strong convergence implies zig-zag spectral convergence]
\label{zz:thm:ZZ-spectral-convergence}
Let $\pi_n:\Gamma\to\Sym(V_n)$ be a sequence of random actions such that $\pi_n |_{\mathds{1}^\perp}$ converge strongly, in probability, to the left regular representation $\lambda$ of $\Gamma$.  Then the following converges in probability
\begin{equation*}
\dH\left(\spec_{\textrm{new}}(G_{\pi_n} \zgzg H)
, \spec(G_\lambda \zgzg H)\right)
\underset{n\to\infty}{\pconv}0.
\label{zz:eq:spectral-general}
\end{equation*}
Consequently, the edges of the new spectra converge.
\end{theorem}

To prove the theorem, we show that $\mA_{G_\lambda \zgzg H}$ and $\mA_{G_{\pi_n} \zgzg H}$ are given by evaluating the same matrix-valued $*$-polynomial, once the generators of $\Gamma$ are represented respectively by $\pi_n$ and by $\lambda$, allowing us to apply strong convergence.
This follows from the following expression of the rotation map in terms of the group action.

\begin{claim}\label{claim:rot_map_formula}
    As a unitary, the rotation map $\mdotG_{\boldsymbol{\pi}}: \ell^2(S) \otimes \ell^2(V) \to \ell^2(S) \otimes \ell^2(V)$ can be written as:
    $$\mdotG_{\boldsymbol{\pi}} = \sum_{s \in S} E_{s,s^{-1}} \otimes \pi(s).$$
\end{claim}

The claim follows directly from comparing the two sides on standard basis vectors.

\begin{proof}[Proof of \cref{zz:thm:ZZ-spectral-convergence}]
Using the expression for the rotation map~(\cref{zz:eq:zz-matrix}), we have
\begin{equation}
\mA _{G_\pi\zgzg H}
=
\sum_{s\in S}
\mA _HE_{s,s^{-1}}\mA _H\otimes\pi(s),
\label{zz:eq:zz-action-form}
\end{equation}
for any action $\pi: \Gamma \to \Sym(V)$. 
Further, if $V$ is finite, then restricting to the orthogonal complement of the constants on the second tensor factor gives the restriction to $\Hzero$ in $\ell^2(S)\otimes \ell^2(V)$:
\begin{equation}
\mA _{G_\pi\zgzg H} |_{\Hzero}
=
\sum_{s\in S}
\mA _HE_{s,s^{-1}}\mA _H\otimes(\pi(s) |_{\mathds{1}^\perp}),
\end{equation}
Enumerate $S = \{s_1, \dots, s_r\}$, and consider the following matrix valued polynomial:
$$ P = \sum_{i=1}^r (\mA_H  E_{s_i, s_i^{-1}} \mA_H) \otimes s_i \in \Mat_{|S|}(\C) \otimes  \C\langle\boldsymbol{x}, \boldsymbol{x}^* \rangle $$
The desired convergence of Hausdorff distances is then a direct application of \cref{lem:matrix_valued_conv} and the assumption that $\pi_n |_{\mathds{1}^\perp}$ converge strongly to $\lambda$.
\end{proof}

Using the proof above, we can also derive the effective bounds in the case of the permutation model.
Notably, the bound holds for an \emph{arbitrary} gadget graph $H$, and considers only its size and degree.

\begin{proposition}[Effective spectral convergence for the permutation model]\label{zz:thm:ZZ-effective-spectral-convergence}
    	Let $\Gamma=\bF_r$ with $r\ge 2$, and let $H$ be a $c$-regular graph on $d=2r$ vertices. For $\Gn$ sampled from the $r$-independent permutations model, we have for every $\epsilon>0$:
    \begin{align*}
        \Prob_{\Gn} \left[ \lambdanew(\Gn \zgzg H) \geq \left\|\mA_{G_\lambda \zgzg H}\right\|_{\op} + \epsilon \right] &\lesssim \frac{d}{n} \left( \frac{d c^2 \log(d/2)}{\epsilon} \right)^8 \log\left( \frac{e \cdot d c^2}{\epsilon} \right)\\
        &= O\left(1/(n\cdot \epsilon^{9})\right).
    \end{align*}
\end{proposition}

\begin{proof}
    This is a direct application of \cref{thm:strong_conv_effective} and the proof of \cref{zz:thm:ZZ-spectral-convergence}. Indeed, if $s_1, \dots, s_r$ are the (non-symmetric) generators of $\bF_r$, in this case the matrix valued polynomial
    $$ P = \sum_{i=1}^r (\mA_H  E_{s_i, s_i^{-1}} \mA_H)  \otimes s_i + \sum_{i=1}^r (\mA_H  E_{s_i^{-1}, s_i} \mA_H)  \otimes s_i^{-1} \in \Mat_{2r}(\C) \otimes  \C\langle\boldsymbol{x}, \boldsymbol{x}^* \rangle $$
is a \emph{linear} matrix polynomial with coefficients in $\Mat_{2r}(\C)$ and we have $\| P \|_1 \leq 2rc^2 $, because:
$$ \| \mA_H  E_{s_i, s_i^{-1}} \mA_H \|_{\op} \leq \|\mA_H \|_{\op} \|  E_{s_i, s_i^{-1}} \|_{\op} \| \mA_H \|_{\op} \leq c^2 $$
where we used that $E_{s_i,s_i^{-1}}$ is a contraction, and that $\| \mA_H \|_{\op} = c$.
\end{proof}

\subsection{The trivial spectrum}\label{zz:sec:trivial-spectrum}
The trivial eigenvalues are easy to compute, and are present in every zig-zag product with a fixed gadget graph $H$, as follows.

\begin{proposition}\label{zz:prop:trivial-block}
For every action $\pi:\Gamma\to\Sym(V)$ on a finite set $V$,
the restriction of $\mA _{G_\pi \zgzg H}$ to ${\Hone}$ is naturally identified with
\begin{equation*}
    \mA _H\mF _S\mA _H
    \qquad\text{on }\ell^2(S).
\end{equation*}
Consequently the spectrum contributed by ${\Hone}$ is
\[
    \spec(\mA _H\mF _S\mA _H).
\]
\end{proposition}

\begin{proof}
We have $\pi(s)\mathds{1}=\mathds{1}$ for every $s\in S$. Hence the restriction of the rotation map to ${\Hone}$ is
\[
\sum_{s\in S}E_{s,s^{-1}}=\mF _S.
\]
Applying \cref{zz:eq:zz-action-form} gives
$\mA _H\mF _S\mA _H$.
\end{proof}

For $\Gamma=\Z_2^{*r}$ one has $\mF _S=I$, so the trivial spectrum is simply
\[
\spec(\mA _H^2).
\]
For $\Gamma=\bF_r$ with its standard symmetric generating set, one recovers the port-flip expression
\[
\spec(\mA _H\mF _S\mA _H).
\]

The following example shows that the port pairing can change both the signs and the magnitudes of trivial eigenvalues.



\begin{example}\label{zz:ex:C6-trivial}
Let $H=C_6$, and number its six ports cyclically by $1,\ldots,6$, identifying
port $i$ with the generator $s_i$. Here $r=3$, so $\bar i=7-i$, and the
port-flip involution and its matrix are
\[
(1\ 6)(2\ 5)(3\ 4),
\qquad
\mF_S=
\begin{pmatrix}
0&0&0&0&0&1\\
0&0&0&0&1&0\\
0&0&0&1&0&0\\
0&0&1&0&0&0\\
0&1&0&0&0&0\\
1&0&0&0&0&0
\end{pmatrix}.
\]
Using the usual cyclic adjacency matrix $\mA_{H}$, direct multiplication gives
\[
\mA_{H}\mF_S\mA_{H}=
\begin{pmatrix}
0&1&0&1&0&2\\
1&0&1&0&2&0\\
0&1&0&2&0&1\\
1&0&2&0&1&0\\
0&2&0&1&0&1\\
2&0&1&0&1&0
\end{pmatrix}.
\]
Direct calculation gives
\[
\spec(\mA_{H}^2)=\{4,4,1,1,1,1\},
\]
and
\[
\spec(\mA_{H}\mF_S\mA_{H})=\{4,1,1,-1,-1,-4\}.
\]

Notice that the trivial spectrum may influence expansion: while the nontrivial spectrum positive edge converges to $\LambdaH=\Lambda_{C_6}=4\cos\left(\frac{\pi}{12}\right)\approx 3.864$ for \emph{both} models, the overall spectrum may or may not have $\lambda_2=4$ as an additional eigenvalue depending on the random model chosen.
\end{example}


The remaining ingredient to proving the main theorems is the following fact.

\begin{proposition}\label{zz:thm:cover-spectrum}
If $H$ is as in \cref{zz:thm:main}, then $\LambdaH$ computes the spectral radius of $\treed \zgzg H$:
$$ \LambdaH = \| \mA_{\treed\zgzg H} \|_{\op} .$$
\end{proposition}

The spectral edge follows from the calculation in \cite[Theorem~5.2 and its proof]{CCM24}. We give the details in \cref{zz:app:cover-spectrum}.

\begin{proof}[Proof of \cref{zz:thm:main} and \cref{zz:thm:main_full_spec}]
\cref{zz:thm:main_full_spec} follows immediately by combining \cref{zz:sec:spectral-convergence} and  \cref{zz:prop:trivial-block}, together with the known strong convergence for the permutation, matching and involution models  \cref{thm:strong_conv_qualitative}.
Together with \cref{zz:thm:cover-spectrum} and the fact that Hausdorff distance between spectra controls the difference between spectral edges, we obtain \cref{zz:thm:main} except for the quantitative assertion in the first clause, which is a direct consequence of \cref{zz:thm:ZZ-effective-spectral-convergence}.

\end{proof}

\begin{remark}[The configuration model]
\label{zz:rem:configuration-model}
    Our zig-zag spectral-convergence results do not cover the configuration model used in the existence proof of \cite{CCM24}. In that model, a single uniformly random perfect matching pairs all $nd$ ports, rather than sampling independent matchings or permutations acting on the $n$ vertices. In general, this randomness does not arise from an action $\pi_n:\Gamma\to\Sym(n)$ and does not admit the fixed matrix-coefficient representation in \cref{zz:eq:zz-action-form} used in our strong-convergence argument.

    Classical contiguity results compare the underlying random regular graphs after the matching decomposition or port-labeling data are forgotten (see \cite{Janson95} and \cite[Section~4.3]{Wormald99}). Since the zig-zag product depends on the rotation map, these comparisons do not directly transfer our spectral-convergence results to the configuration model. Nevertheless, the numerical evidence discussed in \cite[Section~5.4]{CCM24} motivates the following conjecture.
\end{remark}

\begin{conjecture}\label{conj:zigzig_configuration}
	Fix $H$ as above, and let $\Gn$ be sampled from the $d$-regular configuration model, with its induced half-edge labeling. Then, as $n\rightarrow\infty$,
	\[
	\lambda_2(\Gn\zgzg H)\pconv\LambdaH.
	\]
\end{conjecture}

\section{Irregular configuration graphs}\label{sec:configuration}
For each $n$, fix a positive degree sequence
\[
\degseq=(d_1,\ldots,d_n)
\]
with even total degree
\[
N=N_n:=\sum_{i=1}^n d_i.
\]
Assume that the degrees are uniformly bounded: there is a fixed $\maxdeg<\infty$ such that $d_i\le \maxdeg$ for every $n$ and $i$. We also assume that the empirical degree distribution converges to a probability law $\degreeLaw$ on the positive integers. Equivalently, if $U_n$ is a uniformly chosen vertex, then $d_{U_n}$ converges in distribution to $\degreeLaw$. Write
\[
\avgdeg:=\E \degreeLaw,
\qquad
\mindeg:=\min\{d:\Prob(\degreeLaw=d)>0\}.
\]
For example, if $d_1=d_2=\ldots=d_{\frac n2}=3$ and $d_{\frac n2+1}=\ldots=d_n=6$, then $\avgdeg=4.5$ (in this case, for every admissible $n$, $\avgdegfinite=4.5$) and $\mindeg=3$.
Throughout the draft we assume $\mindeg\ge2$. Set $\avgdegn=N_n/n$, so $\avgdegn\longrightarrow\avgdeg$.
We use the standard configuration model \cite{bollobas1980probabilistic,Wormald99}: To sample a graph, attach $d_i$ half-edges to vertex $i$, and choose a uniformly random perfect matching of all $N$ half-edges. Loops and multiple edges are retained.

Let $\Gn$ be the resulting multigraph and $\mAGn$ its adjacency matrix, with each loop contributing $2$ to the diagonal. Thus its row sums are the prescribed degrees. Write
\[
\mDGn:=\diag(d_1,\ldots,d_n).
\]
We study the symmetric normalized adjacency matrix
\[
\mWGn:=\mDGn^{-1/2}\mAGn\mDGn^{-1/2},
\]
which is the appropriate object to study for expansion of irregular graphs. Its eigenvalues are real; write $1=\normeig_1 \ge \normeig_2 \ge \cdots \ge \normeig_n\ge-1$. Smaller values of $\normeig_2$ mean better one-sided spectral expansion. The matrix $\mWGn$ is similar to the stochastic transition matrix $\mPGn=\mDGn^{-1}\mAGn$ of simple random walk. Hence
\begin{equation}
	\frac1n\Tr \mWGn^k
	=\frac1n\Tr \mPGn^k,
	\label{cfg:eq:trace-return}
\end{equation}
which is the average $k$-step return probability from a uniformly chosen vertex.

\subsection{The free-probability benchmark}
\label{cfg:subsec:rotation-free-benchmark}
We recall the following existence result of Cohen and Maor~\cite{CohenMaor2026}, already stated as \cref{thm:cohen-maor-existence-intro} in the introduction.
\begin{theorem}[Existence at the average-degree bound {\cite[Theorem~1]{CohenMaor2026}}]
	\label{cfg:thm:cohen-maor-existence}
	Let $(d_1,\ldots,d_n)$ be a positive degree sequence with even total degree and average $\avgdegfinite\ge2$. There exists a multigraph $G$ in the support of the configuration model with these degrees whose second normalized eigenvalue satisfies
	\begin{equation}\label{eq:cm26}
		\normeig_2(\mWG)\le \frac{2\sqrt{\avgdegfinite-1}}{\avgdegfinite}.
	\end{equation}
\end{theorem}

The fact that $\avgdegfinite$ is the value appearing in the above result can be interpreted in the following way, using the rotation-matrix representation from \cite[Section~5.2]{CohenMaor2026}. Let $\portset_i$ be the set of half-edges attached to vertex $i$, let $\one_{\portset_i}$ be its indicator vector, and define $\mC\in\R^{N\times n}$ by
\[
 \mC e_i=d_i^{-1/2}\one_{\portset_i}.
\]
Let $\mdotGn$ be the rotation matrix of $\Gn$, namely the permutation matrix of the involution that pairs the half-edges. Then
\[
 \mWGn=\mC^\top \mdotGn\mC.
\]
Moreover,
\[
 \mC\mC^\top
 =\diag\left(\frac1{d_1}\mJ_{d_1},\ldots,\frac1{d_n}\mJ_{d_n}\right)
\]
is an orthogonal projection of rank $n$ in dimension $N$, where $\mJ_d$ denotes the $d\times d$ all-ones matrix. Its empirical eigenvalue distribution therefore depends only on $\frac nN = \frac1\avgdegfinite$, not on the detailed degree distribution.

As noted in \cite{CohenMaor2026}, the value in \cref{cfg:thm:cohen-maor-existence} resembles the Ramanujan threshold, so one may ask whether it plays the role of an Alon--Boppana type bound for irregular graphs. However, this is not the case, as \cite{CohenMaor2026} also gives explicit degree sequences for which better expansion is possible.

What remains is the question of typical behavior. The free-probability calculation in \cite{CohenMaor2026} compresses a symmetric $\{\pm1\}$ variable (encoding the eigenvalues of $\mdotGn$) by a free projection of trace $1/\avgdeg$. We denote the resulting benchmark law by $\KM_\avgdeg$. For integer $\avgdeg=d$, this is the normalized Kesten--McKay law \cite{Kesten1959,McKay1981}. Its positive support edge is
$\edgefn(\avgdeg):=2\sqrt{\avgdeg-1}/\avgdeg$.
The family is extended to real $\avgdeg\ge2$ through non-integer free convolution powers (see \cite[Lecture~14]{NS06}).

\subsection{The random limiting law}
\label{cfg:sec:correct}
We next describe the random tree seen from a uniformly chosen vertex, which we think of as the root. Its degree is distributed as $\degreeLaw$. If we traverse an edge, the vertex reached is sampled through one of its half-edges, so higher-degree vertices are more likely to be encountered. Its degree therefore follows the \emph{size-biased} law $\degreeLaw^*$ defined by
\begin{equation}
 \Prob(\degreeLaw^*=c)=\frac{c\Prob(\degreeLaw=c)}{\avgdeg}.
 \label{cfg:eq:size-biased}
\end{equation}
Every nonroot vertex of degree $\degreeLaw^*$ has $\degreeLaw^*-1$ new forward neighbors, and the forward branches are generated independently by the same rule. We denote this random rooted tree by $\Tlim$. This random tree is usually called the unimodular Galton--Watson tree associated with $\degreeLaw$.

The next lemma says that a fixed-size neighborhood in a large random graph is described by $\Tlim$: not only for one random vertex, but also in the proportion of vertices with each possible neighborhood. Here a \emph{radius-$R$ neighborhood} is the subgraph induced by the vertices at distance at most $R$ from the chosen root. This is the bounded-degree case of \cite[Theorem~4.1]{Hofstad2024}; we include the proof for completeness.

\begin{lemma}[Local convergence]
\label{cfg:lem:local-convergence}
Fix a radius $R$ and a finite rooted graph $H$. The proportion of vertices of $\Gn$ whose radius-$R$ neighborhood is isomorphic to $H$ converges in probability to the probability that the radius-$R$ neighborhood of $\Tlim$ is isomorphic to $H$.
\end{lemma}

\begin{proof}
Exploring a fixed-radius neighborhood exposes only $O_{R,\maxdeg}(1)$ vertices and half-edges. At each pairing step, only $O_{R,\maxdeg}(1)$ of the $N=\Theta(n)$ possible partners belong to previously seen vertices. Thus the probability of any collision is $O(1/n)$, and the neighborhood is a tree with probability tending to one.

The root degree converges to $\degreeLaw$ because the root is uniform among vertices. A new vertex is instead reached through a uniform remaining half-edge. If $n_c$ vertices have degree $c$, they initially own $c n_c$ half-edges, and only $O_{R,\maxdeg}(1)$ half-edges have been affected by the exploration. Consequently the probability that the new vertex has degree $c$ is
\[
 \frac{c n_c+O_{R,\maxdeg}(1)}{N+O_{R,\maxdeg}(1)}
 \underset{n\to\infty}{\longrightarrow}\frac{c\Prob(\degreeLaw=c)}{\avgdeg}.
\]
Repeating this calculation gives independent size-biased degrees in the limit, exactly as in $\Tlim$.

It remains to pass from one random root to the proportion of all roots. Let $Y_n$ be the proportion in the statement, and let $p$ be the probability of the corresponding radius-$R$ neighborhood in $\Tlim$. We have already shown that one uniformly chosen root has the desired neighborhood with probability tending to $p$, so $\E Y_n\longrightarrow p$.

Our goal is to show that $Y_n$ itself is close to $p$ with high probability. It is enough to prove that $\operatorname{Var}(Y_n)\longrightarrow0$. Choose two vertices independently and uniformly from $\Gn$. Conditional on $\Gn$, the probability that both have the desired neighborhood is $Y_n^2$, so the unconditional probability is $\E Y_n^2$.
Since degrees are uniformly bounded, the two explorations intersect with probability $O(1/n)$. Conditional on the first exploration, only $O_{R,\maxdeg}(1)$ half-edges have been exposed, and the remaining matching is uniform. Consequently, the second exploration has the same limiting distribution as an independent exploration, regardless of the outcome of the first. Therefore $\E Y_n^2\longrightarrow p^2$.
Since also $\E Y_n\longrightarrow p$,
\[
 \operatorname{Var}(Y_n)=\E Y_n^2-(\E Y_n)^2\longrightarrow0.
\]
Chebyshev's inequality now gives $Y_n\pconv p$, as required.
\end{proof}

\subsubsection{Weak convergence of the empirical spectral measure}
\label{cfg:subsec:limiting-law}
For a rooted graph $(H,r)$, let
\[
 p_k(H,r):=\Prob[X_k=r\mid X_0=r]
\]
be the $k$-step return probability of simple random walk started at $r$. Thus $X_0=r$, each $X_{i+1}$ is chosen uniformly among the neighbors of $X_i$, and we ask for the probability that the walk is back at $r$ at time $k$.

Local convergence also determines the limiting spectral measure. The general principle is developed in \cite[Theorem~1(ii)]{BordenaveLelarge2010} for adjacency and Laplacian matrices.

\begin{proposition}[Limiting spectral law]
\label{cfg:prop:limiting-law}
There is a unique probability measure $\nuD$ on $[-1,1]$ whose moments are
\begin{equation}
 m_k(\nuD):=\int x^k\,d\nuD(x)=\E\left[p_k(\Tlim,r)\right].
 \label{cfg:eq:limiting-moments}
\end{equation}
The empirical spectral measure of $\mWGn$ converges weakly in probability to $\nuD$.
\end{proposition}

\begin{proof}
By \cref{cfg:eq:trace-return}, the $k$-th empirical spectral moment is the average $k$-step return probability over the vertices of $\Gn$. This return probability depends only on the radius-$k$ rooted neighborhood, so \cref{cfg:lem:local-convergence} gives
\[
 \frac1n\Tr \mWGn{}^k
 \underset{n\to\infty}{\pconv}
 \E[p_k(\Tlim,r)]
\]
for every fixed $k$. Since these moments are bounded, the same convergence holds for their expectations. Compactness applied to the expected empirical measures gives a probability measure $\nuD$ with these limiting moments, and moments determine it uniquely on $[-1,1]$. Polynomial approximation then upgrades the displayed moment convergence to weak convergence in probability.
\end{proof}

\subsubsection{The support edge of the limiting spectral law}\label{cfg:sec:support-edge}
The full limiting law depends on the degree distribution, but its support edge has a simple form. For a probability measure $\nu$ on $[-1,1]$, define
\[
\rho(\nu):=\max\{|x|:x\in\supp(\nu)\}.
\]

The needed tree calculation is the following lemma by M{\"u}ller and Wiegel. Applied to the offspring law \(\degreeLaw^*-1\), whose minimum is \(\mindeg-1\), it gives the value
\[
\frac{2\sqrt{\mindeg-1}}{\mindeg}.
\]
The different degree law at the distinguished root of \(\Tlim\) does not change the exponential growth rate of return probabilities, and hence does not change the support edge of the averaged spectral measure \(\nuD\).

\begin{lemma}[{Support edge, \cite[Lemma~6.5]{MuellerWiegel2020}}]\label{cfg:thm:support-edge}
Assume $\mindeg\ge2$. The limiting law $\nuD$ is symmetric and
\begin{equation}
\rho(\nuD)=\rhoDelta:=\frac{2\sqrt{\mindeg-1}}{\mindeg}.
\label{cfg:eq:support-edge}
\end{equation}
Equivalently, $\supp(\nuD)\subseteq[-\rhoDelta,\rhoDelta]$ and both endpoints belong to the support.
\end{lemma}

The point to keep in mind is that the edge is controlled by rare minimum-degree regions, not by a typical neighbor. This is why the answer depends on $\mindeg$ even though the size-biased law $\degreeLaw^*$ favors larger degrees. In contrast, McKenzie \cite{McK26} proves an Alon--Boppana-type lower bound of $\sqrt{\E[\degreeLaw^*-1]}$ for the non-backtracking spectrum of graphs with the same unimodular Galton--Watson local limit. These results are complementary: the non-backtracking bound concerns a different operator and does not directly yield a bound on the normalized random-walk spectrum studied here.

Weak convergence actually gives a statement that does not involve the Cohen--Maor benchmark at all: every fixed threshold below the limiting support edge is exceeded by many eigenvalues of a typical graph.

\begin{corollary}[Thresholds below the limiting edge]
	\label{cfg:cor:atypical}
	For every fixed $t<\rhoDelta$, there exists a constant $c_t>0$ such that, with probability tending to one, at least $c_t n$ eigenvalues of $\mWGn$ are larger than $t$. In particular,
	\[
	\Prob\left(\normeig_2(\mWGn)\le t\right)
	\underset{n\to\infty}{\longrightarrow}0.
	\]
\end{corollary}

\begin{proof}
	Because $t<\rhoDelta$ and $\rhoDelta$ belongs to the support of $\nuD$, there is a continuous function $f\ge0$, vanishing on $[-1,t]$, such that
	\[
	a:=\int f\,d\nuD>0.
	\]
	By weak convergence,
	\[
	\frac1n\sum_{i=1}^n f(\normeig_i)
	\underset{n\to\infty}{\pconv} a.
	\]
	Let $M=\|f\|_\infty$. With probability tending to one, the left-hand side is at least $a/2$. Since only eigenvalues larger than $t$ contribute to the sum, at least $c_t n$ eigenvalues are larger than $t$, where $c_t:=\frac{a}{2M}$. In particular, for all sufficiently large $n$ there is more than one such eigenvalue with probability tending to one, and hence $\normeig_2(\mWGn)>t$ with probability tending to one.
\end{proof}

When $\mindeg=2$, we have $\rhoDelta=1$, and \cref{cfg:cor:atypical}, together with the bound $\normeig_2(\mWGn)\le1$, implies
\[
\normeig_2(\mWGn)\pconv1.
\]

\paragraph{Comparison with \cite{CohenMaor2026} and observed values.}
Suppose that $\degreeLaw$ is nonconstant, so that $\mindeg<\avgdeg$. Since $\edgefn(x)=2\sqrt{x-1}/x$ is strictly decreasing on $[2,\infty)$,
\[
\edgefn(\avgdeg)<\edgefn(\mindeg)=\rhoDelta.
\]
Choose any fixed $t$ strictly between these two values. Since $\avgdegn\to\avgdeg$, for all sufficiently large $n$ the existence theorem \cref{cfg:thm:cohen-maor-existence} guarantees a graph with
\[
\normeig_2\le \edgefn(\avgdegn)<t,
\]
whereas \cref{cfg:cor:atypical} says that a uniformly random configuration graph has $\normeig_2>t$ with probability tending to one. This is the precise sense in which the Cohen--Maor existence guarantee is atypical for a random pairing. The gap can be substantial; \cref{cfg:tab:support-edges} gives some examples.

These results are complementary. \Cref{cfg:thm:cohen-maor-existence} is a finite-$n$ existence theorem for every degree sequence, while our result describes typical random pairings in a bounded-degree asymptotic regime. For $\mindeg\ge3$, the corollary does not establish convergence of $\normeig_2$ to $\rhoDelta$; it only rules out staying below any fixed threshold $t<\rhoDelta$. At accessible sizes, though, the measured second eigenvalue often lies strictly \emph{between} the existence bound and the limiting support edge. We display these experiments, together with the finite-size mechanisms behind them, in \cref{cfg:sec:experiments}.
\begin{table}[ht]
	\centering
	\caption{Equal-mixture degree laws of mean $6$. The existence value from \cite{CohenMaor2026} depends only on the mean, while the support edge of the random limiting law depends on the minimum degree.}
	\label{cfg:tab:support-edges}
	\begin{tabular}{lccc}
		\toprule
		Degree law & $\mindeg$ & Existence bound & Limiting support edge\\
		\midrule
		$6$-regular & $6$ & $0.7454$ & $0.7454$\\
		$\{5,7\}$  & $5$ & $0.7454$ & $0.8000$\\
		$\{4,8\}$  & $4$ & $0.7454$ & $0.8660$\\
		$\{3,9\}$  & $3$ & $0.7454$ & $0.9428$\\
		$\{2,10\}$ & $2$ & $0.7454$ & $1.0000$\\
		\bottomrule
	\end{tabular}
\end{table}

\paragraph{A fourth-moment witness.}
The free-probability benchmark also fails to describe the typical limiting spectral distribution of irregular configuration graphs. Indeed, although $\nuD$ and $\KM_\avgdeg$ have the same second moment, their fourth moments satisfy
\begin{equation}
    m_4(\nuD)-m_4(\KM_\avgdeg)
    =\frac1\avgdeg\left(
    \E\left[\frac1\degreeLaw\right]-\frac1\avgdeg
    \right)^2.
    \label{cfg:eq:m4-gap}
\end{equation}
The right-hand side is strictly positive whenever $\degreeLaw$ is nonconstant, giving an explicit low-order-moment witness to the discrepancy between the free-probability benchmark and the typical random matrix. We derive the identity in \cref{cfg:app:moments}.

\subsection{Numerical experiments and the finite-size puzzle}
\label{cfg:sec:experiments}
For regular configuration graphs, the nontrivial spectral edge converges to the regular-tree value; see \cite[Section~2.2]{Bor20}. Numerically, this limiting picture is already visible at moderate sizes. After identifying the limiting edge for irregular graphs, one might expect a similarly recognizable picture.

However, the numerical picture here is more subtle. The second eigenvalue can remain far below---and hence spectrally better than---the support edge in \cref{cfg:thm:support-edge}. In the reported irregular examples it generally lies \emph{between} the Cohen--Maor existence value and the limiting edge, with only slow upward movement as $n$ grows. As mentioned before, the support edge is controlled by rare minimum-degree regions, and sufficiently large such regions may not occur in a finite graph. We give a simple explanation for this slow finite-size behavior, separately for minimum degree $2$ and minimum degree at least $3$. Experimental details are collected in \cref{cfg:app:experiments}.

\subsubsection{$\mindeg=2$: slow but visible movement toward the edge}
\label{cfg:subsec:deg2-experiment}
\begin{figure}[!htbp]
    \centering
    \includegraphics[width=0.78\textwidth]{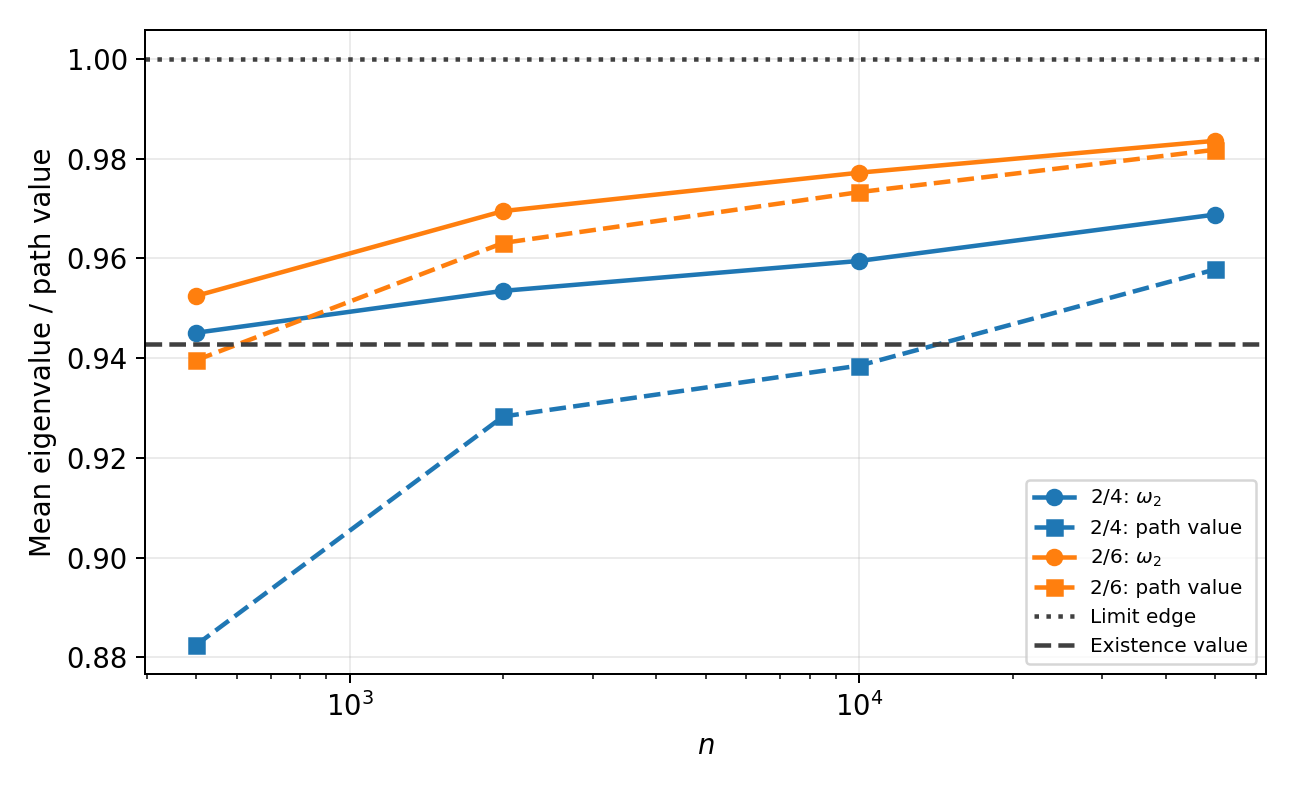}
    \caption{The data of \cref{cfg:tab:deg2-experiment}. Each degree law has one color: solid lines show the average $\normeig_2$ and dashed lines show the average path value from the longest degree-$2$ chain. The horizontal lines mark the limiting edge $1$ and the Cohen--Maor existence value $2\sqrt2/3$.}
    \label{cfg:fig:deg2-convergence}
\end{figure}

We next use two mean-$3$ laws: the $2/4$ family has half the vertices of degree $2$ and half of degree $4$, while the $2/6$ family has three quarters of degree $2$ and one quarter of degree $6$. Both have $\mindeg=2$, so \cref{cfg:thm:support-edge} gives limiting edge $\rho_{2}=1$. The experiments show slow upward movement toward this value, as seen in the solid lines in \cref{cfg:fig:deg2-convergence}.

For comparison with the eigenvalue, let $\longchain(G)$ be the number of vertices in the longest chain of consecutive degree-$2$ vertices in the component being analyzed, with both ends attached to vertices outside the chain. A chain containing $L$ such vertices has the path value
\begin{equation}
 \pathvalue_L:=\cos\left(\frac{\pi}{L+1}\right).
 \label{cfg:eq:path-value}
\end{equation}
A standard test-vector argument on this path gives $\normeig_2(\mWG{})\ge\pathvalue_L-o(1)$ (we include the details for completeness in \cref{cfg:app:path-vector}). Thus long degree-$2$ chains provide concrete finite-size spectral witnesses for large eigenvalues, helping explain the observed experimental results.

\begin{table}[H]

\centering
\caption{Minimum-degree-$2$ experiment. The entries are averages over the sampled graphs.}
\label{cfg:tab:deg2-experiment}
\begin{tabular}{r|ccc|ccc}
\toprule
& \multicolumn{3}{c|}{$2/4$} & \multicolumn{3}{c}{$2/6$}\\
$n$ & $\normeig_2$ & $\longchain$ & path value & $\normeig_2$ & $\longchain$ & path value\\
\midrule
$500$      & .9451 & 5.8  & .8825 & .9525 & 8.3  & .9396\\
$2{,}000$  & .9535 & 7.4  & .9283 & .9695 & 10.8 & .9631\\
$10{,}000$ & .9595 & 8.1  & .9385 & .9772 & 12.9 & .9733\\
$50{,}000$ & .9688 & 9.9  & .9578 & .9836 & 15.7 & .9818\\
\bottomrule
\end{tabular}
\end{table}

The long chains suggest a simple finite-size scale for this slow movement toward the edge. Let $p_2:=\Prob(\degreeLaw=2)$, and let
\[
 q_2:=\Prob(\degreeLaw^*=2)=\frac{2p_2}{\avgdeg}.
\]
Once a path enters a degree-$2$ vertex through one edge, there is only one forward edge. In the local tree model, its endpoint is again degree $2$ with probability $q_2$. Thus a fixed potential chain contains at least $L$ consecutive degree-$2$ vertices with probability approximately $q_2^L$.

There are $O(n)$ possible starting locations. A union bound therefore suggests that chains of length $L$ should not occur once $nq_2^L\ll1$, giving the threshold
\begin{equation}
 \longchain\approx\frac{\log n}{\log(1/q_2)}.
 \label{cfg:eq:deg2-chain-scale}
\end{equation}
For $L$ below this scale one expects many candidate chains, but a matching lower bound requires controlling overlaps. Thus \cref{cfg:eq:deg2-chain-scale} is used here only as a finite-size estimate, not as a theorem. Combining \cref{cfg:eq:deg2-chain-scale,cfg:eq:path-value} suggests the path-based eigenvalue scale
\[
 1-\Theta\left((\log n)^{-2}\right).
\]
The same path obstruction underlies the optimal expansion lower bound of \cite{CohenMaor2026} for graphs of degrees $\{2,3\}$. The data in \cref{cfg:tab:deg2-experiment} show increasing chain lengths and path values, but the sampled sizes are far too small to distinguish the logarithmic scale in \cref{cfg:eq:deg2-chain-scale} from other slowly varying alternatives.

\subsubsection{$\mindeg\ge3$: much slower movement}
\label{cfg:subsec:deg3-experiment}
We now fix mean degree $6$ and minimum degree $3$, while changing how frequently degree $3$ occurs. We use
\[
 \begin{array}{lll}
 3/7:  & \Prob(\degreeLaw=3)=1/4, & q:=\Prob(\degreeLaw^*=3)=0.125,\\
 3/12: & \Prob(\degreeLaw=3)=2/3, & q=0.333,\\
 3/30: & \Prob(\degreeLaw=3)=8/9, & q=0.444.
 \end{array}
\]
Every family has limiting support edge $2\sqrt2/3\approx0.9428$, while the Cohen--Maor existence value at mean $6$ is $2\sqrt5/6\approx0.7454$.

Our experiments reveal values that are neither close to $\rhoDelta$ nor uniformly close to the existence value. They move upward with $n$, and larger $q$ places the finite graph closer to $\rhoDelta$, but the movement is slow.

\begin{table}[H]

\centering
\caption{Median measured $\normeig_2$ for minimum degree $3$. \Cref{cfg:fig:min3-convergence} plots the same experiment, including additional intermediate sizes.}
\label{cfg:tab:deg3-experiment}
\begin{tabular}{r|ccc}
\toprule
$n$ & $3/7$ ($q=.125$) & $3/12$ ($q=.333$) & $3/30$ ($q=.444$)\\
\midrule
$500$      & .7444 & .7646 & .8040\\
$2{,}000$  & .7461 & .7801 & .8086\\
$10{,}000$ & .7550 & .7980 & .8194\\
$50{,}000$ & .7598 & .8094 & .8272\\
\bottomrule
\end{tabular}
\end{table}

\begin{figure}[!htbp]
    \centering
    \includegraphics[width=0.78\textwidth]{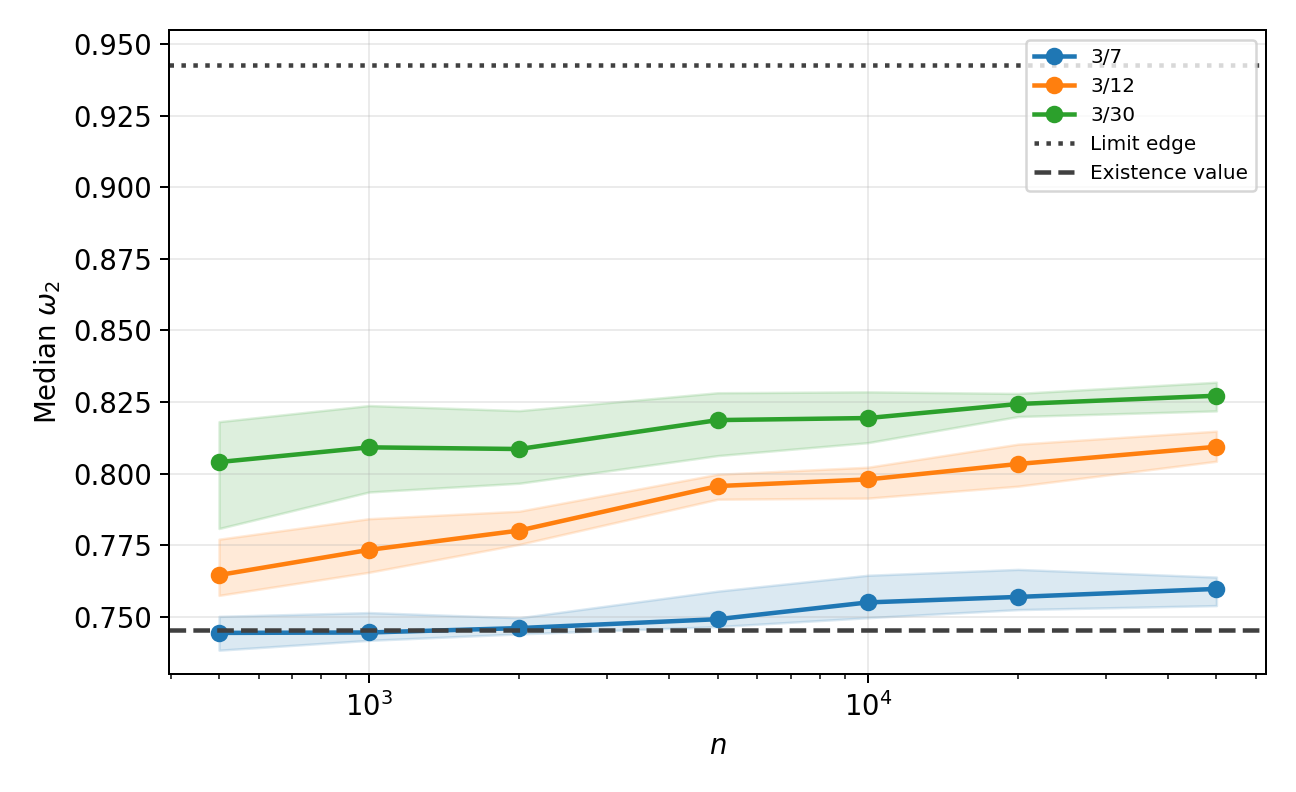}
    \caption{The experiment of \cref{cfg:tab:deg3-experiment} with additional intermediate sizes. Lines show medians and shaded bands show the interquartile range. All three degree laws have mean $6$ and minimum degree $3$, hence the same limiting support edge.}
    \label{cfg:fig:min3-convergence}
\end{figure}

\paragraph{Why is the movement so slow?}
Assume $\degreeLaw$ is nonconstant, and once again set
$q:=\Prob(\degreeLaw^*=\mindeg)\in(0,1)$.
The support edge in \cref{cfg:thm:support-edge} comes from a neighborhood that looks like the $\mindeg$-regular tree; however, in this case this is not a path. After entering such a region through its parent edge, a forward tree with $s$ levels contains
\[
 M_s=1+(\mindeg-1)+\cdots+(\mindeg-1)^{s-1}
 =\frac{(\mindeg-1)^s-1}{\mindeg-2}
\]
vertices whose degrees must all equal $\mindeg$. At a fixed candidate location, the branching model therefore assigns probability $q^{M_s}$ to this event.

Similarly to the $\mindeg=2$ analysis, with $O(n)$ candidate locations, the union-bound threshold $nq^{M_s}\approx1$ yields
\begin{equation}
 M_s\approx\frac{\log n}{\log(1/q)}.
 \label{cfg:eq:Ms-scale}
\end{equation}
Since $M_s$ grows exponentially with $s$,
\begin{equation}
 s_{\max}
 \approx \frac{\log\log n}{\log(\mindeg-1)}+O_{\mindeg,q}(1).
 \label{cfg:eq:loglog-scale}
\end{equation}

As in classical Alon--Boppana arguments, these rare tree regions provide spectral witnesses for large eigenvalues. Consider a forward tree with $s$ levels, all of whose vertices have degree $\mindeg$ in the original graph.
Choose a test vector that assigns to each vertex at distance $j$ from the root the value $j(s-j)(\mindeg-1)^{-j/2}$, for $0\le j<s$. Accounting for the number of vertices at each level, the tree edges alone give a lower bound on the Rayleigh quotient of
\[
\rhoDelta\left(1-\frac{5}{s^2+1}\right)
=\rhoDelta-O_{\mindeg}(s^{-2}).
\]
As in the degree-$2$ case, extending the vector by zero and removing its Perron component introduces only an $O(M_s/n)$ correction. Combining this with the heuristic depth $s_{\max}$ from \cref{cfg:eq:loglog-scale} suggests
\[
\normeig_2\gtrsim
\rhoDelta-O\left((\log\log n)^{-2}\right).
\]

The distinction from degree-$2$ chains is that a minimum-degree tree grows exponentially with its depth, so large regions are much rarer. This provides a simple mechanism for the particularly slow upward movement observed in \cref{cfg:fig:min3-convergence}.

\paragraph{Convergence of the second eigenvalue.}
Our experiments, together with the spectral-witness heuristic, suggest that this slow upward movement continues toward $\rhoDelta$. We therefore formulate the following conjecture.
\begin{conjecture}\label{conj:irregular}
Under the standing assumptions of \cref{sec:configuration}, suppose that $\mindeg\ge3$ and that every degree in every finite sequence is at least $\mindeg$. Then
\[
\normeig_2(\mWGn)\pconv\rhoDelta.
\]
\end{conjecture}
A related Friedman-type question for configuration graphs with fixed degree distributions was raised in \cite{jardon2025skeletons}, in terms of almost-Ramanujan behavior.
A further open problem is to establish rigorous quantitative bounds on the finite-size behavior of $\normeig_2$, and to determine whether the scales suggested by our spectral-witness heuristics are sharp.

\section*{Acknowledgments}
This work was initiated during lecture series given by the first author at the Weizmann Institute of Science and in the Simons Institute for the Theory of Computing in Berkeley. We thank both institutions for their hospitality. We also thank Vedat Levi Alev for fruitful discussions on strong convergence and its applications to expander constructions.
\bibliographystyle{alpha}
\bibliography{references}

\clearpage
\appendix
\crefalias{section}{appendix}
\section{The infinite-cover calculation}\label{zz:app:cover-spectrum}
This appendix records the infinite-cover calculation from \cite[Theorem~5.2 and its proof]{CCM24} in the form used here. We first identify the common limiting graph and then relate first-return walks, all closed walks, and the spectral radius. For a bounded self-adjoint operator $T$, write
\[
    \lamplus(T):=\sup\spec(T),
    \quad
    \lamminus(T):=\inf\spec(T),
\]

\begin{claim}[The common infinite graph]\label{zz:clm:tree-label-independence}
Suppose that $H$ is vertex-transitive. For any two port labelings of the infinite $d$-regular tree $\treed$, the corresponding zig-zag products $\treed\zgzg H$ are isomorphic as graphs. Write
\[
\XH:=\treed\zgzg H.
\]
For every connected port-labeled $d$-regular graph $G$, the graph $\XH$ is a cover of $G\zgzg H$. Moreover, $\XH$ is bipartite, so its spectrum is symmetric about zero. 
\end{claim}

\begin{proof}
Root two port-labeled copies of $\treed$ and proceed recursively away from the root. When a new tree vertex is reached through port $j$ and its image through port $j'$, choose an automorphism of $H$ sending $j$ to $j'$. Vertex-transitivity makes this possible. Map the remaining incident edges according to this automorphism and continue along the tree. Each vertex has a unique parent, so these choices do not conflict. The resulting map identifies the two zig-zag products.
Let $G$ be a connected port-labeled $d$-regular graph, and let $p:\treed\to G$ be its universal covering map. Equip $\treed$ with the port labeling lifted from $G$. Since $p$ preserves the labeled rotation maps, the map
\[
    (s,v)\longmapsto(s,p(v))
\]
induces a covering map from $\treed\zgzg H$ to $G\zgzg H$. By the first part of the claim, the graph on the left is isomorphic to $\XH$.

In the matching and involution models, the regular Schreier graph is $\mathcal T_r$, while in the permutation model it is $\mathcal T_{2r}$. In either case its degree is $d=|V(H)|$. Thus, the corresponding limiting adjacency operators considered in \cref{zz:sec:spectral-convergence} are unitarily equivalent to $\mAXH$, by \cref{zz:eq:zz-matrix} and the preceding isomorphism.
Every zig-zag edge contains exactly one step on $\treed$ and crosses its bipartition.
\end{proof}

For the rest of this appendix, assume that $H$ is as in \cref{zz:thm:cover-spectrum}. This places $H$ within the scope of \cite[Theorem~5.2]{CCM24}. Its proof analyzes a rooted graph $X$, which is a truncation of $\XH$, obtained by deleting finitely many edges and taking the infinite component containing the root $o$.  We use the following consequence of the symbolic and analytic analysis in Sections~5.1.1--5.1.3 there.

Let $\sm{S}_X$ be the class of closed walks in $X$ starting at $o$ and returning to $o$ only at the end, and write
\[
\gf{S}_X(z)=\sum_{\ell\geq1}s_\ell z^\ell,
\]
where $s_\ell$ is the number of such walks of length $\ell$. Similarly, let $\sm{C}_X$ denote the class of all closed walks from $o$, including the empty walk.  Every such walk decomposes uniquely into a sequence of first-return walks, so in the symbolic method notation used in \cite{CCM24},
\begin{equation}
\sm{C}_X=\operatorname{SEQ}(\sm{S}_X),
\qquad
\gf{C}_X(z)=\frac{1}{1-\gf{S}_X(z)},
\label{zz:eq:CX-SX}
\end{equation}
where $\gf{C}_X(z)$ is the analogous generating function.

\begin{lemma}[First-return growth {\cite[Section~5.1]{CCM24}}]\label{zz:lem:TBFTZZP-first-return}
$
\limsup_{\ell\to\infty}s_\ell^{1/\ell}=\LambdaH.
$
\end{lemma}

The equality of the two radii of convergence is stated in \cite[Section~5.1.1]{CCM24}, with a reference to subcriticality.  We give the short verification.

\begin{lemma}[First returns and all returns]\label{zz:lem:first-all-radius}
The generating functions $\gf{S}_X(z)$ and $\gf{C}_X(z)$ have the same radius of convergence.  Consequently, if
\[
\gf{C}_X(z)=\sum_{\ell\geq0}c_\ell z^\ell,
\]
where $c_\ell$ is the number of closed walks of length $\ell$ from $o$, then
\[
\limsup_{\ell\to\infty}c_\ell^{1/\ell}=\LambdaH.
\]
\end{lemma}

\begin{proof}
In \cite[Section~5.1.3]{CCM24}, $\gf{S}_X(z)=E(z^2)$, where $E$ belongs to the smooth inverse-function schema with characteristic value $\tau\in(0,1/c^2)$; see \cite[Corollary~5.9 and Claim~5.11]{CCM24}.  If $\rho$ is the radius of $E$, the full smooth inverse-function theorem gives $E(\rho)=\tau$ \cite[Theorem~VI.6]{FS09}.  Hence, for $R=\sqrt{\rho}$, we have $\gf{S}_X(R)=\tau<1$.  This is precisely the \emph{subcritical} case of the sequence construction; cf. \cite[Section~VI.9]{FS09}.  Indeed, for $|z|<R$,
\[
|\gf{S}_X(z)|\leq \gf{S}_X(|z|)<1,
\]
so $\gf{C}_X(z)=(1-\gf{S}_X(z))^{-1}$ is analytic in $|z|<R$.  Thus $R_{\gf{C}_X}\geq R$.  On the other hand, every first-return walk is also a closed walk, so $c_\ell\geq s_\ell$ for every $\ell$, and therefore $R_{\gf{C}_X}\leq R$.  Hence the radii are equal.  Finally, for a power series $\sum_{\ell\geq0}a_\ell z^\ell$ with radius of convergence $R_a$,
\[
\limsup_{\ell\to\infty} a_\ell^{1/\ell}=\frac{1}{R_a}.
\]
Together with \cref{zz:lem:TBFTZZP-first-return}, this gives the claimed growth rate for $c_\ell$.
\end{proof}

For a connected bounded-degree graph, the norm of its adjacency operator equals the exponential growth rate of closed walks from any fixed vertex.  Since $X$ is bipartite, its spectrum is symmetric about $0$, and \cref{zz:lem:first-all-radius} therefore gives
\begin{equation}
\lamplus(\mAX)=\LambdaH.
\label{zz:eq:truncated-edge}
\end{equation}

It remains only to remove the finite truncation.

\begin{lemma}[Finite truncation]\label{zz:lem:finite-truncation}
Let $X$ be the infinite component obtained from $\XH$ by the finite truncation above.  Then
\[
\lamplus(\mAX)=\lamplus(\mAXH).
\]
\end{lemma}

\begin{proof}
The inequality
\[
\lamplus(\mAX)\leq\lamplus(\mAXH)
\]
is immediate because $X$ is a subgraph of $\XH$.

For the reverse inequality, note that $\XH$ is vertex-transitive: the same rooted-tree recursion used in \cref{zz:clm:tree-label-independence}, now started from any prescribed pair of vertices of $\XH$, produces an automorphism carrying one to the other. The graph $\XH$ need not be connected. For example, when $H=C_6$ and we use the matching labeling, every zig-zag step preserves the parity of the port, so $\XH$ is disconnected. Nevertheless, vertex-transitivity implies that all connected components of $\XH$ are isomorphic.
Consequently, $\lamplus(\mAXH)$ can be approximated by finitely supported test vectors within any chosen component, in particular the one containing $X$. Fix $\varepsilon>0$ and choose a finitely supported unit vector $f$ in that component such that
\[
\langle f,\mAXH f\rangle
\geq
\lamplus(\mAXH)-\varepsilon.
\]
The support of $f$ lies in some finite-radius ball.  Since only finitely many edges were deleted and $X$ is infinite, $X$ contains a vertex whose ball of that radius is disjoint from all deleted edges.  By vertex-transitivity, move the support of $f$ to that ball.  The Rayleigh quotient is unchanged, and the translated vector is supported inside $X$.  Hence
\[
\lamplus(\mAX)
\geq
\lamplus(\mAXH)-\varepsilon.
\]
Letting $\varepsilon\to0$ proves the claim.
\end{proof}

\begin{proof}[Proof of \cref{zz:thm:cover-spectrum}]
By \cref{zz:eq:truncated-edge,zz:lem:finite-truncation},
\[
\lamplus(\mAXH)=\LambdaH.
\]
By \cref{zz:clm:tree-label-independence}, the graph $\XH$ is bipartite, so its spectrum is symmetric about zero. Therefore,
\[
    \lamminus(\mAXH)=-\LambdaH,
    \quad
    \|\mAXH\|_{\op}=\LambdaH.
\]
\end{proof}

\section{The degree-2 path test vector}
\label{cfg:app:path-vector}
This appendix proves the finite-path bound used in \cref{cfg:subsec:deg2-experiment}. The long-path obstruction appears in \cite[Lemma~6.4]{MuellerWiegel2020} and \cite[Theorem~5]{CohenMaor2026}. We include the normalized-adjacency calculation and its finite-graph orthogonality correction for completeness.

Suppose $G$ is connected, has $n$ vertices, and contains a chain $v_1,\ldots,v_L$ of $L=o(n)$ consecutive degree-$2$ vertices, with the two ends attached to vertices outside the chain. The principal normalized-adjacency matrix on these $L$ vertices is one half of the path adjacency matrix, so its largest eigenvalue is
\[
 \pathvalue_L=\cos\left(\frac{\pi}{L+1}\right).
\]
Let $x$ be its unit top eigenvector, extended by zero outside the chain, and let $p$ be the unit Perron vector of $G$. Set $\alpha=\langle x,p\rangle$. Since $p_v=\sqrt{d_v/N}$,
\[
 \alpha^2\le\sum_{i=1}^L p_{v_i}^2=\frac{2L}{N}.
\]
The vector $y=x-\alpha p$ is orthogonal to $p$, and
\[
 \frac{y^\top \mWG{}y}{\|y\|^2}
 =\frac{\pathvalue_L-\alpha^2}{1-\alpha^2}
 =\pathvalue_L-O(L/n)
\]
for bounded degrees. Hence $\normeig_2(\mWG{})\ge\pathvalue_L-O(L/n)$ whenever $L=o(n)$.

\section{A fourth-moment comparison}
\label{cfg:app:moments}

This appendix derives the fourth-moment identity in \cref{cfg:eq:m4-gap} by enumerating short closed walks on the limiting tree $\Tlim$.

The free-probability benchmark $\KM_\avgdeg$ has first two even moments
\begin{equation}
    m_2(\KM_\avgdeg)=\frac1\avgdeg,
    \qquad
    m_4(\KM_\avgdeg)=\frac{2\avgdeg-1}{\avgdeg^3}.
    \label{cfg:eq:free-moments}
\end{equation}
On the limiting tree $T=\Tlim$, an edge $u\sim v$ has normalized weight
\[
(\mWT)_{uv}=\frac1{\sqrt{d_ud_v}}.
\]
For a closed walk $v_0\to v_1\to\cdots\to v_k=v_0$, its contribution to $(\mWT^k)_{v_0v_0}$ is
\[
\prod_{i=0}^{k-1}\frac1{\sqrt{d_{v_i}d_{v_{i+1}}}}
=\prod_{i=0}^{k-1}\frac1{d_{v_i}}.
\]
The equality holds for every closed walk: each occurrence of a degree appears once through the incoming edge and once through the outgoing edge, so the square roots pair up. Summing these weights gives the return probability to $v_0$; averaging over the law of the rooted tree $(\Tlim,r)$ gives $m_k(\nuD)$.

The second moment agrees with the free-probability benchmark. Conditional on a root of degree $c$, there are $c$ choices for the intermediate neighbor, and the factor $1/c$ at the root cancels when we sum over them. Thus, by \cref{cfg:eq:size-biased},
\[
m_2(\nuD)=\E \left[\frac1{\degreeLaw^*}\right]=\sum_c \frac1c \cdot \frac{c\Prob(\degreeLaw=c)}{\avgdeg}=\frac1\avgdeg=m_2(\KM_\avgdeg).
\]
Following a random edge favors larger degrees, but for a two-step return the factor $1/\degreeLaw^*$ cancels this bias. Repeated visits, however, introduce higher inverse powers of the encountered degrees, where this cancellation no longer occurs. The fourth moment is already enough to see the difference.

\paragraph{The fourth moment.}
Write
\[
\hinv:=\E\left[\frac1\degreeLaw\right].
\]
We will use two elementary consequences of the size-biased law. From \cref{cfg:eq:size-biased},
\begin{align}
	\E\left[\frac1{\degreeLaw^*}\right]
	&=\sum_c \frac1c\frac{c\Prob(\degreeLaw=c)}{\avgdeg}
	=\frac1\avgdeg,
	\label{cfg:eq:size-biased-first}\\
	\E\left[\frac1{(\degreeLaw^*)^2}\right]
	&=\sum_c \frac1{c^2}\frac{c\Prob(\degreeLaw=c)}{\avgdeg}
	=\frac1\avgdeg\sum_c\frac{\Prob(\degreeLaw=c)}c
	=\frac \hinv\avgdeg.
	\label{cfg:eq:size-biased-second}
\end{align}
Because $\Tlim$ is a tree, a closed walk of length four has only two shapes.

\paragraph{Two backtracks from the root.}
The walk is $v\to u\to v\to w\to v$. Condition on the root degree $\degreeLaw=c$, and let $X_i=1/d_{u_i}$ for the $c$ neighbors of the root. The total contribution of all ordered choices $(u,w)$ is
\[
\frac1{c^2}\left(\sum_{i=1}^c X_i\right)^2.
\]
The $c$ terms with the same neighbor use $\E X_i^2=\hinv/\avgdeg$; the $c(c-1)$ terms with different neighbors use $\E X_iX_j=\E X_i \cdot\E X_j=1/\avgdeg^2$. Thus the conditional expectation is
\[
\frac1c\frac \hinv\avgdeg+\left(1-\frac1c\right)\frac1{\avgdeg^2}.
\]
Averaging over $\degreeLaw$ gives
\begin{equation}
	A=\frac{\hinv^2}{\avgdeg}+\frac{1-\hinv}{\avgdeg^2}.
	\label{cfg:eq:typeA}
\end{equation}

\paragraph{A depth-two excursion.}
The walk is $v\to u\to x\to u\to v$ with $x\ne v$. If $u$ has degree $c$, there are $c-1$ possible forward choices for $x$. Averaging over $\degreeLaw^*$ gives
\begin{align*}
	B
	&=\frac1\avgdeg\E\left[\frac{\degreeLaw^*-1}{(\degreeLaw^*)^2}\right]\\
	&=\frac1\avgdeg\left(
	\E\!\left[\frac1{\degreeLaw^*}\right]
	-\E\!\left[\frac1{(\degreeLaw^*)^2}\right]\right)
	=\frac{1-\hinv}{\avgdeg^2}.
\end{align*}
Overall,
\begin{equation}
	m_4(\nuD)=A+B=\frac{\hinv^2}{\avgdeg}+\frac{2(1-\hinv)}{\avgdeg^2}.
	\label{cfg:eq:m4-conf}
\end{equation}
Subtracting the free-probability fourth moment in \cref{cfg:eq:free-moments} yields
\[
    m_4(\nuD)-m_4(\KM_\avgdeg)
    =\frac1\avgdeg\left(\hinv-\frac1\avgdeg\right)^2,
\]
which is precisely \cref{cfg:eq:m4-gap}, since $\hinv=\E[1/\degreeLaw]$. The expression is nonnegative. Moreover, strict convexity of $x\mapsto1/x$ gives
\[
    \hinv=\E\left[\frac1\degreeLaw\right]
    >\frac1{\E[\degreeLaw]}=\frac1\avgdeg
\]
whenever $\degreeLaw$ is nonconstant. Thus the fourth-moment discrepancy is strictly positive for every nonconstant degree distribution.

\section{Experimental details and numerical checks}
\label{cfg:app:experiments}
The multigraph configuration model was sampled by uniformly matching half-edges; loops and multiple edges were retained. For the extreme-eigenvalue experiments we used the largest connected component, so that extra eigenvalue-$1$ multiplicities from small disconnected components do not obscure the nontrivial edge. The experimental questions and comparisons were chosen by the authors, with ChatGPT assisting in the implementation and execution of the numerical experiments.

In the minimum-degree-$2$ experiment of \cref{cfg:tab:deg2-experiment,cfg:fig:deg2-convergence}, the displayed values of $\normeig_2$, $\longchain$, and the path value are averages over the sampled graphs. The path-value column averages $\pathvalue_{\longchain}$ separately for each graph, rather than applying $\pathvalue$ to the averaged value of $\longchain$. In the minimum-degree-$3$ experiment of \cref{cfg:tab:deg3-experiment,cfg:fig:min3-convergence}, the lines show medians and the shaded bands show the interquartile range.

As a numerical check of the fourth-moment calculation in \cref{cfg:app:moments}, \cref{cfg:tab:m4-experiment} compares the fourth moment at $n=100{,}000$ with the limiting value in \cref{cfg:eq:m4-conf}. Weak convergence itself is proved analytically in \cref{cfg:prop:limiting-law}; the table is only a numerical check. In contrast with the slow finite-size behavior of the extreme eigenvalues in \cref{cfg:subsec:deg2-experiment,cfg:subsec:deg3-experiment}, the fourth moment is already close to its limiting prediction at this size.

\begin{table}[ht]
\centering
\caption{Fourth-moment check at $n=100{,}000$.}
\label{cfg:tab:m4-experiment}
\begin{tabular}{l|cc|cc}
\toprule
& \multicolumn{2}{c|}{Analytic} & \multicolumn{2}{c}{Experiment}\\
Family & $m_4(\KM_\avgdeg)$ & $m_4(\nuD)$ & mean & median\\
\midrule
$3$-regular & .18519 & .18519 & .18519 & .18519\\
$2/4$       & .18519 & .18576 & .18577 & .18576\\
$6$-regular & .05093 & .05093 & .05094 & .05094\\
$4/8$       & .05093 & .05100 & .05101 & .05102\\
\bottomrule
\end{tabular}
\end{table}

\section*{AI methodology}
The proofs of the random zig-zag product results in \cref{sec:main-zigzag} were developed entirely without AI assistance. For the results on irregular graphs in \cref{sec:configuration} and its supporting appendices, we used OpenAI's ChatGPT to assist with developing proofs, locating relevant previous results, and implementing and executing numerical experiments. The experimental questions and comparisons were chosen by the authors.

\end{document}